\input amstex
\documentstyle{amsppt}
\magnification=1200 \hsize=13.8cm \catcode`\@=11
\def\NoLogo{\let\logo@\empty}
\catcode`\@=\active \NoLogo

\def\heat{\lf(\frac{\p}{\p t}-\Delta\ri)}
\def \b {\beta}

\def\lf{\left}
\def\ri{\right}
\def\bbar{\bar \beta}
\def\a{\alpha}

\def\g{\gamma}
\def\e{\epsilon}
\def\p{\partial}
\def\delbar{{\bar\delta}}
\def\ddbar{\partial\bar\partial}
\def\dbar{\bar\partial}

\def\C{\Bbb C}
\def\R{\Bbb R}

\def\vp{\varphi}

\def\dbar{\bar\partial}

\def\abb{{\alpha\bar\beta}}
\def\gbd{{\gamma\bar\delta}}

\def \D {\Delta}
\def\aint{\frac{\ \ }{\ \ }{\hskip -0.4cm}\int}
\documentstyle{amsppt}
\magnification=1200 \hsize=13.8cm \vsize=19 cm

\leftheadtext{Lei Ni} \rightheadtext{A monotonicity formula}
\topmatter
\title{A monotonicity formula on complete K\"ahler manifolds with nonnegative bisectional
curvature}\endtitle

\author{Lei Ni\footnotemark }\endauthor
\footnotetext"$^{1}$"{Research partially supported by NSF grant
DMS-0328624, USA.}
\address
Department of Mathematics, University of California, San Diego, La
Jolla, CA 92093
\endaddress
\email{ lni\@math.ucsd.edu}
\endemail

\affil { University of California, San Diego}
\endaffil

\date  July 2003\enddate

\abstract In this paper, we derive a new monotonicity formula for
the plurisuhbarmonic functions on complete K\"ahler manifolds with
nonnegative bisectional curvature. As applications we derive the
sharp estimates for the dimension of the spaces of holomorphic
functions (sections) with polynomial growth, which in particular,
partially solves a conjecture of Yau.
\endabstract

\endtopmatter

\document

\subheading{\S0 Introduction}\vskip .2cm

In \cite{Y}, Yau proposed to study the uniformization of complete
K\"ahler manifolds with nonnegative curvature. In particular, one
wishes to determine  if a complete K\"ahler manifold $M$ with
positive bisectional curvature is biholomorphic to $\C^m$ or not.
For this sake, it was further asked in \cite{Y} that if the the
ring of the holomorphic functions with polynomial growth, which we
denote by ${\Cal O}_P(M)$, is finitely generated or not,  and if
the dimension of the spaces of holomorphic functions of polynomial
growth is bounded from above by the dimension of the corresponding
spaces of polynomials on $\C^m$ or not. This paper addresses the
later question. We denote by ${\Cal O}_d(M) $ the space of
holomorphic functions of polynomial growth with degree $d$. (See
Section 3 for the precise definition.) Then ${\Cal
O}_P(M)=\bigcup_{d\ge 0}{\Cal O}_d(M)$. In this paper, we show
that

 \proclaim {Theorem 0.1} Let $M^m$ be a complete K\"ahler
manifold with nonnegative holomorphic bisectional curvature.
Assume that $M$ is of maximum volume growth. Then
$$
\dim_{\C}({\Cal O}_d(M))\le \dim_{\C}({\Cal O}_{[d]}(\C^m)). \tag
0.1
$$
Here $[d]$ is the greatest integer less than or equal to $d$.
\endproclaim

Denote $V_o(r)$ the volume of the ball of radius $r$ centered at
$o$. For manifolds with nonnegative Ricci curvature,
$\frac{V_o(r)}{r^{2m}}$ is monotone decreasing by the Bishop
volume comparison theorem. $M$ is called to have the maximum
volume growth if $\lim_{r\to \infty} \frac{V_o(r)}{r^{2m}}>0$.

Although we have not quite proved the finite generation for the
ring ${\Cal O}_P(M)$, we can show that the quotient field
generated by ${\Cal O}_P(M)$ is finitely generated. In fact, we
obtain the following dimension estimate for the general case.

\proclaim{Theorem 0.2} Let $M^m$ be a complete K\"ahler manifold
with nonnegative holomorphic bisectional curvature of complex
dimension $m$. There exists a constant $C_1=C_1(m)$ such that for
every $d\ge 1$,
$$
\dim_{\C}({\Cal O}_d(M))\le C_1 d^m. \tag 0.2
$$
\endproclaim

By an argument in \cite{M}, which is originally due to Poincar\'e
and Siegel, the above result does imply that the rational
functions field ${\Cal M}(M)$ generated by ${\Cal O}_P(M)$ is of
transcendental degree at most $m$. From this one can further
construct a birational embedding of $M$ into $\C^{m+2}$ in the
case $M$ has positive bisectional curvature and admits nonconstant
holomorphic functions of polynomial growth. In a future
publication we shall study the finite generation of ${\Cal
O}_P(M)$ as well as the affine embedding of $M$, using the results
and techniques developed here.

The new idea of this paper is a monotonicity formula for the
plurisubharmonic functions (as well as positive currents). In
order to illustrate our approach let us recall a classical result
attributed as  Bishop-Lelong Lemma.

\medskip

 {\it Let $\Theta$ be a $(p,p)$ positive current in $\C^m$.
Define
 $$
\nu(\Theta, x,r)=\frac{1}{r^{2m-2p}}\int_{B_x(r)}\Theta
\wedge\lf(\frac{1}{\pi} \omega_{\C^m}\ri)^{m-p}.\tag 0.3
$$
Here $\omega_{\C^m}$ is the K\"ahler form of $\C^m$. Then
$$
\frac{\p}{\p r} \nu(\Theta, x,r)\ge 0. \tag 0.4
$$
}

\medskip
\noindent  This monotonicity formula in particular can be applied
to the $(1,1)$ current  $\frac{\sqrt{-1}}{2\pi}\ddbar \log |f|^2$,
where $f$ is a holomorphic function.  Through the monotonicity
(0.4), in \cite{B} Bombieri derived a
 Schwartz's Lemma type  inequality, out of which one can infer
 that the vanishing order of a polynomial is bounded by its
 degree.

 However, this line of argument encountered  difficulties,
 when applied to the non-flat spaces.  \cite{M} made the first such attempt.
 The following result of Mok
  in \cite{M} is particularly notable.

\proclaim{Theorem (Mok)} Let $M$ be a complete K\"ahler manifold
with nonnegative holomorphic bisectional curvature. Suppose that
there exists  positive constants $C_2$ and $C_3$ such that for
some fixed point $o\in M$
$$
V_o(r)\ge C_2r^{2m} \tag 0.5
$$
and
$$
0<{\Cal R}(x)<\frac{C_3}{(1+r(x))^2}. \tag 0.6
$$
Here $V_o(r)$ is the volume of $B_o(r)$, the ball of radius $r$
centered at $o$, ${\Cal R}(x)$ is the scalar curvature function
and $r(x)$ is the distance function to $o$. Then $M$ is
biholomorphic to an affine algebraic variety.
\endproclaim

The key step of the proof to the above result is to obtain
estimate (0.2) and a multiplicity estimate, from which (0.2) can
be  derived. (See Section 3 for the definition and derivation.)
The extra assumption (0.5) and (0.6) were needed to compensate the
failure of (0.4) on curved manifolds.

The main contribution of this paper is to establish a new
monotonicity formula on any complete K\"ahler manifolds with
nonnegative bisectional curvature. The monotonicity formula was
established through  the heat equation deformation of the initial
plurisubharmonic functions (or positive $(1,1)$-current). In the
case of plurisubharmonic functions, the monotonicity formula has
the following simple form. One can refer Theorem 1.1 for the
general case.

\medskip

 {\it Let $M$ be a complete K\"ahler manifold. Let
$v(x,t)$ be a family of plurisubharmonic functions deformed by the
heat equation $\heat v(x,t)=0$ such that $w(x,t)=v_t(x,t)$ is
continuous for each $t>0$. Then
$$
\frac{\p}{\p t} \lf(tw(x,t)\ri)\ge 0. \tag 0.7
$$
 }

\medskip

\noindent Here we assume that $v(x,t)$ is plurisubharmonic just
for the sake of simplicity.  This assumption  in general can be
ensured by a recent established maximum principle for tensors on
complete manifolds in \cite {NT2, Theorem 2.1}, if the initial
function $v(x,0)$ is plurisubharmonic and of reasonable growth.
The monotonicity of $tw(x,t)$ replaces (0.4) in the non-flat case.
The dimension estimates in Theorem 0.1 and Theorem 0.2 can be
proved by comparing the value of $tw(x,t)$ at $t=0$ with its limit
 as $ t\to \infty$. In the proof of Theorem 0.1, the sharp upper
bound on the heat kernel by Li-Tam-Wang \cite{LTW, Theorem 2.1}
was used. In the proof of Theorem 0.2, we make use of the general
`moment' estimates proved in \cite{N1, Theorem 3.1} by the author.

The estimate (0.7) follows from  a gradient estimate of Li-Yau
type, which resembles the trace form of Hamilton's Li-Yau-Hamilton
differential inequality's \cite{H}, originally also called the
differential Harnack inequality,  for the Ricci flow. Indeed, the
derivation of (0.7) was motivated by the earlier work of
Chow-Hamilton in \cite {CH} on the linear trace Li-Yau-Hamilton
inequality, as well as \cite{NT1} by
 Luen-Fai Tam and the author for the K\"ahler case. In fact,
 the author discovered (0.7) when
 trying to generalize Chow's  interpolation \cite{C} between the Li-Yau's
 estimates and the linear trace Harnack (Li-Yau-Hamilton) estimates for the
 Ricci flow on Riemann surfaces to the high dimension. It is also influenced by a
 discussion held with G. Perelman, in which Perelman attributed the
 reason for the success of parabolic method to an `uncertainty
 principle'. This might  suggest an elliptic method may  only
 be  possible after deeper understandings of the geometry of
K\"ahler
 manifolds with nonnegative curvature, such as a total
 classification of such manifolds up to biholomorphisms. Since one can think
 the heat equation deformation of a plurisubharmonic functions is a parabolic deformation
 of related currents, the work
 here suggests that there exist strong connections between the
 K\"ahler-Ricci flow and the other curvature flows.
 The recent works of Perelman \cite{P} and Huisken-Sinestrari \cite{HS} also suggest some
 strong dualities between the Ricci flow and the mean curvature
 flow. It is not clear if the parabolic deformation of the
 currents in this paper has any connection with the mean curvature
 flow or not. This certainly deserves further
 deeper investigations in the future projects. The previous work
 \cite{NT2} is also crucial to this paper, especially the tensor
 maximum principle on complete manifolds \cite{Theorem 2.1, NT2}.

There are many works on estimating the dimension of the harmonic
functions of polynomial growth in the last a few years. See, for
example \cite{CM, LT1-2, LW}.  One can refer \cite{L} for a
survey. The previous results on harmonic functions  conclude that
the dimension has upper bound of the form $C_5\, d^{2m-1}$, which
is sharp in the power for the harmonic functions. Since the space
of harmonic functions is far bigger than the space of holomorphic
functions, the estimate is not sharp for the holomorphic
functions. While estimate (0.2) is sharp in the power and strong
enough to draw some complex geometric conclusions out of it.

In \cite{LW}, the problem of obtaining the sharp upper bound for
the dimension of the space of  harmonic functions of polynomial
growth was studied for manifolds with nonnegative sectional
curvature and maximum volume growth. An asymptotically sharp
estimate was proved there. But the estimate as (0.1) is still
missing. Due to the apparent difference of the nature of the two
problems  the method in this paper is quite different from the
previous papers on harmonic functions. (The exceptional cases are
either $m=1$ or $d=1$.  For both cases, the sharp bounds have been
proved by Li-Tam \cite{LT1-2}.)

Combining the estimate (0.3) and the H\"ormander's $L^2$-estimate
of $\dbar$-operator we can obtain some topological consequences on
the complete K\"ahler manifolds with nonnegative bisectional
curvature. For example, we have the following result.

\proclaim{Theorem 0.3} Let $M$ be a complete K\"ahler manifold
with nonnegative bisectional curvature. Assume that the
transcendence degree of ${\Cal M}(M)$
$$
deg_{tr}({\Cal M}(M))=m. \tag 0.8
$$
Then $M$ has finite fundamental group.
\endproclaim

 Since it is still unknown if a complete K\"ahler manifold with
 positive bisectional curvature is simply-connected or not, the
 result above gives some information towards this  question as well as the uniformization problem.
 The assumption (0.8) can be replaced by the positivity of the Ricci
 and some average curvature decay conditions.

We organize the paper as follows. In Section 1 we derive the
gradient estimate for the symmetric tensors, from which the
monotonicity formula (0.7) is derived in Section 2. Theorem 0.1
and Theorem 0.2 are proved in Section 4 and Section 3
respectively.  The more general  estimates are also  proved for
the holomorphic sections of polynomial growth for line bundles
with finite `Lelong number at infinity'. (See Section 4 for the
precise definition.) As another application, we also include
unified treatments on  the Liouville theorem for the
plurisubharmonic functions on complete K\"ahler manifolds with
nonnegative bisectional curvature (namely any continuous
plurisubharmonic functions of $o(\log r)$ growth is a constant),
as well as  the optimal gap theorem \cite{NT2, Corollary 6.1}.
They can all be phrased as the positivity of the `Lelong number at
infinity' for  non-flat, nonnegative holomorphic line bundles.
These results, which are presented in  Section 2 as a warm-up to
the later cases, were originally proved in \cite{NT2} by Luen-Fai
Tam and the author, using different methods. In fact, Corollary
6.1 in \cite{NT2} is more general than what proved in Theorem 2.2
here. In Section 5 we prove the finite generation of the quotient
field and some topological consequences.

\medskip

{\it Acknowledgment.} The author would like to thank Professors
Peter Li, Luen-Fai Tam and  Jiaping Wang for helpful discussions
and their interests to this work. The author would also like to
thank Professors Ben Chow,  John Milnor and the Clay foundation
for supporting his trip in April to Simons' lectures at SUNY Stony
Brook.

\input amstex
\documentstyle{amsppt}
\magnification=1200 \hsize=13.8cm \catcode`\@=11
\def\NoLogo{\let\logo@\empty}
\catcode`\@=\active \NoLogo
\def\heatt{\lf (\Delta-\frac{\p}{\p t}\ri)}

\def\heat{\lf(\frac{\p}{\p t}-\Delta\ri)}

\def \b {\beta}

\def\lf{\left}
\def\ri{\right}
\def\bbar{\bar \beta}
\def\a{\alpha}

\def\g{\gamma}
\def\e{\epsilon}
\def\p{\partial}
\def\delbar{{\bar\delta}}
\def\ddbar{\partial\bar\partial}
\def\dbar{\bar\partial}

\def\C{\Bbb C}
\def\R{\Bbb R}

\def\vp{\varphi}

\def\dbar{\bar\partial}

\def\abb{{\alpha\bar\beta}}
\def\gbd{{\gamma\bar\delta}}

\def \D {\Delta}
\def\aint{\frac{\ \ }{\ \ }{\hskip -0.4cm}\int}
\vsize=19.0 cm

\subheading{\S1 A gradient estimate}

\vskip .2cm

In this section we derive a new gradient estimate for the
symmetric (1,1) tensor $h_{\abb}(x,t)$ satisfying the
Lichnerowicz-Laplacian heat equation:
$$
\heat h_{\gbd}=R_{\beta\bar{\a}\gbd
}h_{\abb}-\frac{1}{2}\lf(R_{\gamma
\bar{s}}h_{s\delbar}+R_{s\delbar}h_{\gamma \bar{s}}\ri). \tag 1.1
$$
We assume that $h_{\abb}(x,t)\ge 0$ and $M$ has nonnegative
bisectional curvature. Applying the estimate to the complex
Hessian of a plurisubharmonic function we can obtain a parabolic
version the classical three-circle theorem for the subharmonic
functions on the complex plane. The condition $h_{\abb}(x,t)\ge 0$
can be ensured in most case if $h_{\abb}(x,0)\ge 0$ by the maximum
principle proved recently in \cite{NT2}.

For any (1, 0) vector field $V$ we define
$$
\split
 Z&=\frac{1}{2}\lf(g^{\abb}\nabla_{\bar{\beta}}div
(h)_{\a}+g^{\gbd}\nabla_{\gamma}div(h)_{\delbar}\ri)+g^{\abb}div(h)_{\a}V_{\bar{\beta}}+g^{\gbd}
div(h)_{\delbar}V_{\gamma}\\
&\quad \quad
+g^{\abb}g^{\gbd}h_{\a\delbar}V_{\bar{\beta}}V_{\gamma}+\frac{H}{t}.
\endsplit
\tag 1.2
$$
Here
$$div(h)_\a=g^{\gbd}\nabla_{\gamma}h_{\a\delbar}, \quad
div(h)_{\delbar}=g^{\abb}\nabla_{\bar{\beta}}h_{\a\delbar} \tag
1.3
$$
and
$$
H=g^{\abb}h_{\abb}. \tag 1.4
$$

 \proclaim{Theorem 1.1} Let $M$ be a complete K\"ahler manifold
 with nonnegative holomorphic bisectional curvature. Let
 $h_{\abb}(x,t)\ge 0$ be a symmetric (1,1) tensor satisfying
 (1.1) on $M\times (0, T)$.  Assume that for any $\epsilon'>0$,
 $$
\int_{\e'}^{T}\int_M e^{-ar^2(x)}\|h\|^2\, dv\, dt <\infty. \tag
1.5
$$
Then
$$
Z(x,t)\ge 0, \tag 1.6
$$
for any (1,0) vector $V$.
 \endproclaim

In order to prove the theorem let us first state some simple
lemmas.

\proclaim{Lemma 1.1}
$$
\heat div(h)_\a =-\frac{1}{2}R_{\a\bar{t}}\, div(h)_t, \ \ \ \heat
div(h)_{\bar{\a}}=-\frac{1}{2}R_{\bar{\a}t}\, div(h)_{\bar{t}}.
\tag 1.7
$$

\endproclaim
\demo{Proof} Direct calculation shows that
$$
\split \frac{\p}{\p t}\lf (
g^{\g\delbar}\nabla_{\g}h_{\a\delbar}\ri) & =
g^{\g\delbar}\frac{\p}{\p t}\lf(\p_\g
h_{\a\delbar}-\Gamma^{p}_{\a\g}
h_{p\delbar}\ri)\\
& = \nabla_{\g}\lf(\D h_{\a\bar{\g}}+R_{\a\bar{\g}
s\bar{t}}h_{\bar{s}t} -\frac{1}{2}R_{\a\bar{t}}h_{t\bar{\g}}-
\frac{1}{2}R_{t\bar{\g}}h_{\a\bar{t}}\ri)\\ &= \nabla_\a
R_{s\bar{t}}h_{\bar{s}t} +R_{\a\bar{\g}s\bar{t}}\nabla_\g
h_{\bar{s}t}-
\frac{1}{2}\nabla_{\g}R_{\a\bar{t}}h_{t\bar{\g}}\\
& \ \ -\frac{1}{2}R_{\a\bar{t}}\nabla_{\g}h_{t\bar{\g}}-
\frac{1}{2}\nabla_t Rh_{\a\bar{t}}-
\frac{1}{2}R_{t\bar{\g}}\nabla_{\g}h_{\a\bar{t}}+\nabla_{\g} (\D
h_{\a\bar{\g}}).
\endsplit \tag1.8
$$
Now we calculate $\nabla_{\g}(\D h_{\a\bar{\g}})$. By definition,
$$
\nabla_{\g}(\D h_{\a\bar{\g}}) =\frac{1}{2}\nabla_{\g}
\lf(\nabla_s\nabla_{\bar{s}}+\nabla_{\bar{s}}\nabla_s\ri)h_{\a\bar{\g}}.
$$
On the other hand,
$$
\split \nabla_{\g} \nabla_{s} \nabla_{\bar{s}} h_{\a\bar{\g}} & =
\nabla_s \nabla_\g \nabla_{\bar{s}} h_{\a\bar{\g}} \\
& = \nabla_{s} \lf[ \nabla_{\bar{s}} \nabla_\g  h_{\a\bar{\g}}-
R_{\a\bar{p}\g\bar{s}}h_{p\bar{\g}}+R_{p\bar{\g}\g\bar{s}} h_{\a\bar{p}}\ri]\\
&= \nabla_{s}\nabla_{\bar{s}} \nabla_\g h_{\a\bar{\g}}- \nabla_\g
R_{\a\bar{p}}h_{p\bar{\g}}-R_{\a\bar{p}\g\bar{s}} \nabla_s
h_{p\bar{\g}}+\nabla_p R h_{\a\bar{p}}+ R_{p\bar{s}}\nabla_s
h_{\a\bar{p}}.
\endsplit
$$
Similarly,
$$
\split \nabla_\g \nabla_{\bar{s}}\nabla_s h_{\a\bar{\g}} & =
\nabla_{\bar{s}} \nabla_\g\nabla_s
h_{\a\bar{\g}}+R_{p\bar{\g}\g\bar{s}}
\nabla_s h_{\a\bar{p}}-R_{s\bar{p}\g\bar{s}}\nabla_p h_{\a\bar{\g}}\\
& \ \ \ - R_{\a\bar{p}\g \bar{s}}\nabla_s h_{p\bar{\g}}\\
&= \nabla_{\bar{s}}\nabla_s\nabla_{\g}h_{\a\bar{\g}}+
R_{p\bar{s}}\nabla_s h_{\a\bar{p}}-R_{\g\bar{p}}\nabla_p
h_{\a\bar{\g}} -R_{\a\bar{p}\g \bar{s}}\nabla_s h_{p\bar{\g}}.
\endsplit
$$
Combining the above three equations we have that
$$
\split \nabla_{\g}(\D h_{\a\bar{\g}}) &= \D(\nabla_\g
h_{\a\bar{\g}})-\frac{1}{2} \nabla_\g R_{\a\bar{p}}h_{p\bar{\g}}-
R_{\a\bar{p}\g\bar{s}}\nabla_s h_{p\bar{\g}} \\
& \ \ \ +\frac{1}{2}\nabla_p R h_{\a\bar{p}}+ R_{p\bar{s}}\nabla_s
h_{\a\bar{p}}-\frac{1}{2}R_{\g\bar{p}}\nabla_p h_{\a\bar{\g}}.
\endsplit
$$
Plugging the above equation into (1.8), the first equation in the
lemma is proved. The second one is the conjugation of the first.
\enddemo
\proclaim{Lemma 1.2}
$$
\heat \lf(g^{\abb}\nabla_{\bar{\beta}} div(h)_{\a}\ri)=0,\ \ \ \
\heat \lf(g^{\abb}\nabla_{\beta} div(h)_{\bar{\a}}\ri)=0. \tag 1.9
$$
\endproclaim
\demo{Proof} Follows from Lemma 1.1 and routine calculations. In
deed,
$$
\split \frac{\p}{\p t}\lf(g^{\a\bbar}\nabla_{\bbar}div(h)_\a\ri) &
=\nabla_{\bar{\a}}\lf(\frac{\p}{\p t} div(h)_\a\ri)\\
 & = \nabla_{\bar{\a}}\lf[\D div(h)_{\a}- \frac{1}{2}R_{\a\bar{t}}div(h)_t\ri],
\endsplit
$$
by Lemma 1.1. Therefore we have that
$$
\split \frac{\p}{\p t}\lf(g^{\a\bbar}\nabla_{\bbar}div(h)_\a\ri)
&= \nabla_{\bar{\a}}\lf(\D div(h)_\a\ri)-\frac{1}{2}R_{s\bar{\a}}
\nabla_{\bar{s}}\lf(div(h)_\a\ri)-
\frac{1}{2}\nabla_{\bar{t}}R\lf( div(h)_t\ri).
\endsplit \tag1.10
$$
Now we calculate $\nabla_{\bar{\a}}\lf(\D div(h)_\a\ri)$. By
definition
$$
\nabla_{\bar{\a}}\lf(\D div(h)_\a\ri)=\frac{1}{2}\nabla_{\bar{\a}}
\nabla_s\nabla_{\bar{s}}div(h)_{\a}+\frac{1}{2}\nabla_{\bar{\a}}
\nabla_{\bar{s}}\nabla_s div(h)_\a.
$$
On the other hand
$$
\split \nabla_{\bar{\a}}\nabla_{\bar{s}}\nabla_s div(h)_\a &=
\nabla_{\bar{s}}\nabla_{\bar{\a}}\nabla_s div(h)_\a \\
& = \nabla_{\bar{s}}\lf[ \nabla_s\nabla_{\bar{\a}}div(h)_\a+
R_{\a\bar{p}s\bar{\a}}div(h)_p\ri]\\
&= \nabla_{\bar{s}}\nabla_s\nabla_{\bar{\a}}div(h)_\a+
\lf(\nabla_{\bar{s}}R\ri)\lf(
div(h)_s\ri)+R_{s\bar{p}}\nabla_{\bar{s}}div(h)_p
\endsplit
$$
and
$$
\split \nabla_{\bar{\a}}\nabla_s\nabla_{\bar{s}}div(h)_\a &=
\nabla_{s}\nabla_{\bar{\a}}\nabla_{\bar{s}}div(h)_{\a} +
R_{\bar{p}s}\nabla_{\bar{s}} div(h)_p-R_{p\bar{\a}}\nabla_{\bar{p}}div(h)_\a\\
&= \nabla_{s}\nabla_{\bar{s}}\nabla_{\bar{\a}}div(h)_{\a}+
R_{\bar{p}s}\nabla_{\bar{s}}
div(h)_p-R_{p\bar{\a}}\nabla_{\bar{p}}div(h)_\a.
\endsplit
$$
Combining the above three equalities we have that
$$
\nabla_{\bar{\a}}\lf(\D
div(h)_\a\ri)=\D\lf(\nabla_{\bar{\a}}div(h)_{\a}\ri)
+\frac{1}{2}\nabla_{\bar{s}}R(div(h)_s)+\frac{1}{2}R_{s\bar{p}}
\nabla_{\bar{s}}div(h)_p.
$$
Plugging into (1.10), this completes the proof of the first
equation of Lemma 1.2. The second one is the conjugation of the
first.

\enddemo

Since $h_{\abb}\ge 0$ only (not strictly positive), we use the
perturbation trick as in \cite{NT1}. Namely, we consider
$$
\split
 \widehat Z&=\frac{1}{2}\lf(g^{\abb}\nabla_{\bar{\beta}}div
(h)_{\a}+g^{\gbd}\nabla_{\gamma}div(h)_{\delbar}\ri)+g^{\abb}div(h)_{\a}V_{\bar{\beta}}+g^{\gbd}
div(h)_{\delbar}V_{\gamma}\\
&\quad \quad +g^{\abb}g^{\gbd}\lf(h_{\a\delbar}+\e g_{\a
\delbar}\ri)V_{\bar{\beta}}V_{\gamma}+\frac{H+\e m}{t}.
\endsplit
\tag 1.11
$$
We can simply denote $\widetilde h_{\abb}=h_{\abb}+\e g_{\abb}$,
which is strictly positive definite. Let $V$ be the vector field
which minimizes $\widehat Z$. Then the first variation formula
gives
$$
div(h)_\a +\widetilde h_{\a\bar{\g}}V_\g =0\  \text{  and  }\
div(h)_{\bar{\a}} +\widetilde h_{\g\bar{\a}}V_{\bar{\g}}=0.
\tag1.12
$$
Differentiate (1.12) we have that
$$
\split & \nabla_s div(h)_\a+\lf(\nabla_s h_{\a\bar{\g}}\ri)V_{\g}
+\widetilde h_{\a\bar{\g}}\nabla_s V_\g =0, \ \ \nabla_s
div(h)_{\bar{\a}} +\lf(\nabla_s h_{\g\bar{\a}}\ri)
V_{\bar{\g}}+\widetilde h_{\g\bar{\a}}\nabla_s
V_{\bar{\g}}=0,\\
&
\nabla_{\bar{s}}div(h)_\a+\lf(\nabla_{\bar{s}}h_{\a\bar{\g}}\ri)V_{\g}
+\widetilde h_{\a\bar{\g}}\nabla_{\bar{s}} V_\g =0,\ \
\nabla_{\bar{s}} div(h)_{\bar{\a}} +\lf(\nabla_{\bar{s}}
h_{\g\bar{\a}}\ri) V_{\bar{\g}}+ \widetilde
h_{\g\bar{\a}}\nabla_{\bar{s}} V_{\bar{\g}}=0.
\endsplit \tag1.13
$$
From (1.12) we also have the following alternative form of
$\widehat Z$
$$
\widehat Z=-\frac{1}{2}\widetilde
h_{\a\bbar}\nabla_{\bar{\a}}V_{\beta} -\frac{1}{2}\widetilde
h_{\beta\bar{\a}}\nabla_{\a}V_{\bbar}+\frac{H+\e m}{t}. \tag 1.14
$$

On the other hand, by Lemma 1.1, 1.2 we have that
$$
\split \lf(\frac{\p}{\p t}-\D \ri)\widehat Z & =  div(h)_\a \lf(
\lf(\frac{\p}{\p t}-\D \ri)V_{\bar{\a}}\ri)+div(h)_{\bar{\a}}
\lf(\lf(\frac{\p}{\p t}-\D \ri)V_\a\ri)\\
& \ \ \ -\nabla_s
div(h)_\a\nabla_{\bar{s}}V_{\bar{\a}}-\nabla_{\bar{s}}
div(h)_\a\nabla_s V_{\bar{\a}}\\
&\ \ \ -\nabla_s div(h)_{\bar{\a}}
\nabla_{\bar{s}}V_{\a}-\nabla_{\bar{s}}div(h)_{\bar{\a}}\nabla_s V_\a\\
& \ \ \ - \frac{1}{2}R_{\a\bar{t}}div(h)_t V_{\bar{\a}}
-\frac{1}{2}R_{t\bar{\a}}div(h)_{\bar{t}}V_\a
\\
& \ \ \ +R_{\abb s\bar{t}}h_{\bar{s}t}V_{\beta}V_{\bar{\a}}-
\frac{1}{2}R_{\a\bar{s}}h_{s\bar{\g}}V_{\g}V_{\bar{\a}}
-\frac{1}{2}h_{\a\bar{s}}R_{s\bar{\g}}V_\g V_{\bar{\a}}\\
& \ \ \  + \widetilde h_{\a\bar{\g}}\lf(\lf(\frac{\p}{\p t}-\D
\ri)V_\g\ri)V_{\bar{\a}}
+\widetilde h_{\a\bar{\g}}V_\g \lf(\lf(\frac{\p}{\p t}-\D \ri)V_{\bar{\a}}\ri)\\
& \ \ \ -\nabla_s h_{\a\bar{\g}}\nabla_{\bar{s}}\lf(V_\g
V_{\bar{\a}}\ri)
-\nabla_{\bar{s}}h_{\a\bar{\g}}\nabla_{s}\lf(V_{\g}V_{\bar{\a}}\ri)\\
& \ \ \  - \widetilde h_{\a\bar{\g}}\lf[\nabla_s V_\g
\nabla_{\bar{s}}V_{\bar{\a}}+
\nabla_{\bar{s}}V_{\g}\nabla_{s}V_{\bar{\a}}\ri] -\frac{H+\e
m}{t^2}.
\endsplit
\tag1.15
$$
\proclaim{Lemma 1.3} $$ \split \heat \widehat Z &= \widetilde
h_{\gamma
\bar{\a}}\lf[\nabla_pV_{\bar{\gamma}}-\frac{1}{t}g_{p\bar{\gamma}}\ri]\lf[\nabla_{\bar{p}}V_{\a}
-\frac{1}{t}g_{p\bar{\a}}\ri]+\widetilde h_{\gamma
\bar\a}\nabla_{\bar{p}}V_{\bar{\gamma}}\nabla_pV_{\a}\\
&\ \ +R_{\abb
s\bar{t}}h_{\bar{s}t}V_{\beta}V_{\bar{\a}}-\frac{2\widehat
Z}{t}.\endsplit \tag 1.16
$$
\endproclaim
\demo{Proof} Using (1.12), (1.13) we can simplify (1.15) to
$$
\split \lf(\frac{\p}{\p t}-\D \ri)\widehat Z & = R_{\abb
s\bar{t}}h_{\bar{s}t}V_{\beta}V_{\bar{\a}}\\
&\ \ \ +\widetilde h_{\g\bar{\a}}\nabla_s
V_{\bar{\g}}\nabla_{\bar{s}}V_{\a} +\widetilde
h_{\g\bar{\a}}\nabla_{\bar{s}}V_{\bar{\g}}\nabla_s V_\a
-\frac{H+\e m}{t^2}.
\endsplit
\tag1.17
$$
Combining with (1.14) we have the lemma.
\enddemo

In order to apply the maximum principle on complete manifolds and
prove the theorem we need to the following result.

\proclaim{Lemma 1.4} Under the assumption of Theorem 1.1, for any
$\eta>0$,
$$
\int_{\eta}^{T} \int_M e^{-ar^2(x)}\|div(h)\|^2\, dv\, dt <\infty
\tag 1.18
$$
and
$$
\int_{\eta}^{T} \int_M e^{-ar^2(x)}\lf(\|\nabla_{s}
div(h)_{\a}\|^2+\|\nabla_{\bar{s}}div(h)_{\a}\|^2\ri) dv\, dt
<\infty. \tag 1.19
$$
\endproclaim
\demo{Proof} To simplify the notation we first define
$$
\split & \Phi =\|h\|^2,\\
&\Psi =\|div(h)\|^2,\\
&\Lambda=\|\nabla_s
div(h)_{\a}\|^2+\|\nabla_{\bar{s}}div(h)_{\a}\|^2. \endsplit
$$
From (1.1), we have  that
$$
\heatt \Phi \ge \Psi. \tag 1.20
$$
Here we have used the fact that $M$ has nonnegative holomorphic
bisectional curvature. The reader can refer \cite{NT2, Lemma 2.2}
for a detailed proof of this fact. It follows from an argument
which goes back to Bishop and Goldberg \cite{BG}. (See also
\cite{MSY}.)

Let $\phi$ be  a cut-off function such that $\phi=0$ for $r(x)\ge
2R$ or  $t\le \frac{\eta}{2}$ and $\phi=1$ for $r(x)\le R$ and
$t\ge \eta$. Multiplying $\phi^2$ on both side of (1.20),
integration by parts gives that
$$
\int_0^T\int_M \Psi \phi^2\, dv \, dt\le 2\int_0^T\int_M \Phi
\lf(|(\phi^2)_t|+4|\nabla \phi|^2\ri)\, dv\, dt.
$$
Then (1.18) follows from the assumption (1.5).

To prove (1.19) we need to calculate $\heatt \Psi$. From Lemma
1.1, it is easy to obtain
$$
\split \heatt \Psi &= \Lambda
+R_{\a\bar{t}}div(h)_{\bar{\a}}div(h)_t \\
 & \ge \Lambda
 \endsplit \tag 1.21
$$
since $M$ has nonnegative Ricci curvature. Repeat the above
argument for (1.18). We complete the proof of the lemma.

\enddemo

\demo{Proof of Theorem 1.1} By translating the time and the
limiting argument we can assume that $h_{\abb}$ is well-defined on
$M \times[0, t]$. Since $\widetilde h_\abb\ge \epsilon g_\abb$ on
$M\times[0,T]$, by (1.12) and (1.13), we have
$$
||V||\le C_1\Psi^{\frac{1}{2}},
$$
and
$$
||\nabla V||\le C_2\lf(\Lambda^{\frac{1}{2}}+\Psi\ri),
$$
for some constants $C_1$ and $C_2$. Combining this with (1.14), we
have that
$$
|t^2\widehat Z|^2\le C_3t\lf(\Phi
+\Phi(\Psi^2+\Lambda)+1\ri)\tag1.22
$$
for some constant $C_3$.   By (1.16), the corresponding
$t^2\widehat Z$
 satisfies  that $$
\lf(\frac{\p}{\p t}-\D \ri)(t^2\widehat Z)\ge0 \tag1.23
$$
for the vector field which minimizes $\widehat Z$. By  Lemma 1.4
and (1.22), we have that
$$
\int_0^{T}\int_M\exp(-ar^2_0(x))\lf(t^2\widehat Z\ri)^2\, dv\,
dt<\infty
$$
for any $a>0$. By the maximum principle of Karp-Li \cite{KL} (See
also \cite{NT1, Theorem 1.1}), we have $t^2\widehat Z\ge0$ because
it is obvious that $t^2\widehat Z=0$ at $t=0$. Since this is true
for the vector field $V$ minimizing $\widehat Z$, we have
$\widehat Z\ge0$ for any
 (1,0) vector field. Let $\epsilon\to0$ and the proof of the theorem is completed.
\enddemo

\proclaim{Remark} Theorem 1.1 was motivated by the so-called
linear trace Li-Yau-Hamilton inequality for Ricci flow. The linear
trace Li-Yau-Hamilton inequality of the real case was first proved
in \cite{CH} by Chow and Hamilton. In \cite{NT1}, the authors
proved  the corresponding one for the K\"ahler-Ricci flow. In
fact, we can state Theorem 1.2 in \cite{NT1} in a slight more
general way such that it can be used in classifying  the
K\"ahler-Ricci solitons, which was done in \cite{N2}.
\endproclaim

\proclaim{Theorem 1.2} Let $(M, g_{\abb}(x,t))$ be a complete
solution to the K\"ahler-Ricci flow on $M\times(0, T)$ with
nonnegative bisectional curvature. Assume that the curvature is
uniformly bounded on $M\times\{t\}$ for any $t>0$. Let $h$ be a
symmetric (1,1) tensor satisfying (1.1). Assume also that
$h_{\abb}(x,t)\ge 0$ and (1.5) holds. Then
$$
\widetilde Z \ge 0 \tag 1.24
$$
where
$$
\split \widetilde Z
&=\frac12[g^{\a\bbar}\nabla_{\bbar} div(h)_\a+ g^{\g\delbar}\nabla_\g div(h)_\delbar]\\
&\quad+g^{\a\bbar}g^{\g\delbar}[R_{\a\delbar}h_{\g\bbar}+\nabla_\g
h_{\a\delbar}V_{\bbar}+\nabla_{\bbar}h_{\a\delbar}V_{\g}
+h_{\a\delbar}V_{\bbar}V_\g]+\frac{H}{t}.
\endsplit \tag 1.25
$$
\endproclaim
Besides relaxing the assumption on $h$, another main advantage of
stating result as above is that the form here  is more useful
without considering the initial value problem.  The form  stated
in \cite{CH}, as well as in  \cite{NT1},  with the initial value
prevents the application to the expanding solitons. It is also
more clear to separate the issue of preserving the nonnegativity
of $h$ from the nonnegativity of $\widetilde Z$.

\input amstex
\documentstyle{amsppt}
\magnification=1200 \hsize=13.8cm \catcode`\@=11
\def\NoLogo{\let\logo@\empty}
\catcode`\@=\active \NoLogo
\def\heatt{\lf (\Delta-\frac{\p}{\p t}\ri)}

\def\heat{\lf(\frac{\p}{\p t}-\Delta\ri)}

\def \b {\beta}

\def\lf{\left}
\def\ri{\right}
\def\bbar{\bar \beta}
\def\a{\alpha}

\def\g{\gamma}
\def\e{\epsilon}
\def\p{\partial}
\def\delbar{{\bar\delta}}
\def\ddbar{\partial\bar\partial}
\def\dbar{\bar\partial}

\def\C{\Bbb C}
\def\R{\Bbb R}
\def\tv{\tilde v}
\def\vp{\varphi}

\def\dbar{\bar\partial}

\def\abb{{\alpha\bar\beta}}
\def\gbd{{\gamma\bar\delta}}

\def \D {\Delta}
\def\aint{\frac{\ \ }{\ \ }{\hskip -0.4cm}\int}
\vsize=19.0 cm

\subheading{\S2 Nonnegative holomorphic line bundles}

\vskip .2cm

In this section we shall apply results in Section 1 to study the
holomorphic line bundles on K\"ahler manifolds with nonnegative
holomorphic bisectional curvature. First we illustrate the cases
when Theorem 1.1 can be applied.

\proclaim{Theorem 2.1} Let $(E, H)$ be a holomorphic vector bundle
on $M$. Consider the Hermitian metric $H(x,t)$ deformed by the
Hermitian-Einstein flow:
$$
\frac{\p H}{\p t}H^{-1}=-\Lambda F_H +\lambda I.\tag 2.1
$$
Here $\Lambda$ means the contraction by  the K\"ahler form
$\omega$, $\lambda$ is a constant, which is a holomorphic
invariant in the case $M$ is compact, and $F_H$ is the curvature
of the metric $H$, which locally can be written as
$F^j_{i\abb}dz^{\a}\wedge d\bar{z}^{\b}e_{i}^*\otimes e_j$ with
$\{e_i\}$ a local frame for $E$. The transition rule for $H$ under
the frame change is
$H^U_{i\bar{j}}=f_i^k\overline{f_j^k}H^V_{k\bar{l}}$ with
transition functions $f_i^j$ satisfying $e^U_i=f_i^je^V_j$. Denote
$\rho=\frac{\sqrt{-1}}{2\pi}\Omega_{\abb}dz_{\a}\wedge
d\bar{z}_{\b}=\frac{\sqrt{-1}}{2\pi}\sum_{i}F^i_{i\abb}
dz^{\a}\wedge d\bar{z}^{\b}$. Assume that $\Omega_{\abb}$ is
smooth on $M\times(0, t]$. Then $\Omega_{\abb}(x,t)$ satisfies
(1.1). Therefore, if $\Omega_{\abb}(x,t)\ge 0$, (1.6) holds under
the assumption of (1.5). In particular, if $\Omega(x, t)>0$,
$$
\Omega_t-\frac{|\nabla \Omega|^2}{\Omega}+\frac{\Omega}{t}\ge 0.
\tag 2.2
$$
\endproclaim
\demo{Proof} See pages 10--12 of \cite{N2} for the proof.
\enddemo

The following Harnack inequality follows from an argument of
Li-Yau \cite{LY} and (2.2).
 \proclaim{Corollary 2.1} Let $M$ and $\Omega$ be as above.  Assume
 that $\Omega(x,t)>0$. Then for any $t_2>t_1$
 $$
\Omega(x,t_2)\ge \Omega(y,
t_1)\lf(\frac{t_1}{t_2}\ri)\exp\lf(-\frac{r^2(x,
y)}{4(t_2-t_1)}\ri). \tag 2.3
 $$
 In particular,
 $$
\frac{\p}{\p t}\lf(t\, \Omega(x,t)\ri) \ge 0. \tag 2.4
 $$
\endproclaim

The Hermitian-Einstein flow (2.1) was studied, for example by
Donaldson \cite{Dn}, to flow a metric into an equilibrium
 solution under some algebraic stability assumptions. In this
 section we focus on the following two cases when Theorem 2.1
 applies.

{\it Case 1:} In the special case $E=L$ is a line bundle and
$\lambda =0$, (2.1) reduces to the simple equation:
$$
\heat v (x,t)=\Omega(x). \tag 2.5
$$
and $\Omega_{\abb}(x,t)=\Omega_{\abb}(x)+v_{\abb}(x,t)$.  Here
$h(x,t)=h(x)\exp(-v(x,t))$ with $v(x,0)=0$ solves (2.1) with
$h(x,0)=h(x)$. It is easy to see that
$w(x,t)=v_t(x,t)=\Omega(x,t)$ satisfies the heat equation $\heatt
w(x,t)=0$ with the initial data $w(x,0)=\Omega(x)$. In the
following we will focus ourself to the line bundle case.

We call $(L, h)$ is nonnegative if the curvature of $(L, h)$,
$\rho=\frac{\sqrt{-1}}{2\pi}\Omega_{\abb}\, dz_{\a}\wedge
d\bar{z}_\beta =\frac{\sqrt{-1}}{2\pi}\ddbar \log (h)$, is
nonnegative $(1,1)$ form.

In order to ensure that $\Omega_{\abb}(x,t)\ge 0$ when
$\Omega_{\abb}(x,0)=\Omega_{\abb}(x)\ge 0$ we need to put some
constraints on $\Omega(x)$. First we assume $\Omega(x)$ is
continuous.  Furthermore we also require  that
$$
\sup_{r\ge 0}\lf(\exp(-ar)\aint_{B_o(r)}\Omega(y)\, dy\ri)<\infty
\tag 2.6
$$
and
$$
\sup_{r\ge 0}\lf( \exp(-ar^2)\aint_{B_o(r)}\Omega^2(y)\,
dy\ri)<\infty \tag 2.7
 $$
for some positive constant $a>0$. Here $\aint_{B_o(r)}$ means the
$\frac{1}{V_o(r)}\int_{B_o(r)}$.

{\it Case 2:} The  case when $\Omega_{\abb}(x,0)$ is given by the
Hessian of  continuous plurisubharmonic functions has special
interest. This corresponds the singular metrics
$h(x)=\exp(-u(x))$, as those considered in \cite{D}, since we do
not require the smoothness on $u(x)$. (In \cite{D}, dual to the
local nature, the locally integrable functions are allowed.)
However since we are interested in global properties and our
argument is global we require the functions to be defined on whole
$M$. Moreover in order to apply the tensor maximum principle
\cite{Theorem 2.1, NT2} we also put growth constrains on the
plurisubharmonic functions instead of (2.6)and (2.7), which are
specified as follows.

Let $u$ be a continuous function on $M$. We call $u$ is of
exponential growth if for some $a>0 $ such that
$$
|u|(x)\le \exp(a (r^2(x)+1)). \tag 2.6'
$$
By Proposition 2.1 of \cite{NT1} we know that if $u(x)$ is of
exponential growth $\heatt \tv=0$ with $\tv(x,0)=u(x)$ has
solution on $M\times[0, T]$ for any $T>0$. Furthermore, we know
that there exists  a constant $b$ such that
$$
|\tv|(x,t)\le \exp(b (r^2(x)+1)). \tag 2.8
$$
In this case, it is easy to see that $h(x,t)=\exp(-\tv(x,t))$
gives the solution to (2.1) and $v(x,t)=\tv(x,t)-u(x)$ solves
(2.5) with $\Omega_{\abb}(x,t)=\tv_{\abb}(x,t)$.

The following lemma ensures $\Omega_{\abb}(x,t)\ge 0$ for the
above two cases.

\proclaim{Lemma 2.1} Let $M$ be a complete K\"ahler manifold with
nonnegative holomorphic bisectional curvature. let $(L, h)$ is a
nonnegative holomorphic line bundle. We assume that either we are
in case 1 with (2.6) and (2.7), or in case 2 with (2.6'). Then
(2.1) has long time solution with $\Omega_{\abb}(x,t)\ge 0$.
\endproclaim
\demo{Proof} The {\it Case 2} is easier. Since (2.1) amounts to
solving $\heatt \tv=0$ with $\tv(x,0)=u(x)$, the result follows
from Theorem 3.1 of \cite{NT2}.

For the {\it Case 1}, clearly,
$$
v(x,t)=\int_0^t \int_M H(x,y, s) \Omega(y)\, dv_y \, ds
$$
gives the solution to (2.5). It exists for all time due to (2.6).
In order to show $\Omega_{\abb}(x,t)\ge 0$, since
$\Omega_{\abb}(x,t)$ satisfies (1.1) by Theorem 2.1, we only need
to check that the maximum principle \cite{NT2, Theorem 2.1}
applies. Due to the assumption (2.7), we only need to verify it
for $v_{\abb}$. Notice that $v(x,t)$ which satisfies the
non-homogeneous heat equation (2.5). Therefore $v(x,t)$ has
point-wise control through the representation formula above. The
standard integration by parts arguments give the wanted integral
estimates for $\|v_{\abb}\|^2$. The interested reader can refer
the proof of Lemma 6.2 of \cite{NT2} for details of checking on
conditions for the maximum principle. See also Lemma 1.4. The
extra terms caused by the  non-homogeneous term $\Omega(x)$ will
be taken care of by the assumption (2.7).
\enddemo

Combining Theorem 2.1 and Lemma 2.1 we are in the position to
apply the monotonicity formula (2.4). In the next  we prove the
following gap theorem, which combines the Liouville theorem
\cite{NT2, Theorem 0.3} with  a special case of \cite{Corollary
6.1, NT2}. The proof here uses the result from Section 1 and fits
the general duality principle in \cite{N2}.

\proclaim{Theorem 2.2} Let $M$ be a complete K\"ahler manifold
with nonnegative holomorphic bisectional curvature. Let $(L, h)$
be a holomorphic line bundle on $M$ with hermitian metric $h$. We
assume either in the Case 1 that (2.7) holds and
$$
\int_0^r s\lf(\aint_{B_o(s)}\Omega(y)\, dy\ri)\, ds =o(\log r),
\tag 2.9
$$
for some $a>0$, with $\Omega(y)=g^{\abb}\Omega_{\abb}(y)$; or in
the Case 2,
$$
\limsup_{r\to \infty}\frac{u(x)}{\log r} =0. \tag 2.10
$$
Then $(L, h)$ is flat. Namely $\Omega_{\abb}\equiv 0$. In
particular, if $L=K_M^{-1}$, the anti-canonical line bundle, it
implies that $M$ is flat. Moreover, in the Case 2 it further
implies $u$ is a constant.
\endproclaim

\demo{Proof of Theorem 2.2} Clearly in the first case (2.9)
implies (2.6) and  in the second case, since we can replace $u$ by
$u_{+}$, the positive part of $u$, we have (2.6'). Therefore we
can apply Theorem 2.1, in particular (2.4) in both situations, by
Lemma 2.1.

Assume that $(L, h)$ is not flat. Then $\Omega(x)\ge 0$ and $>0$
somewhere. This implies that $\Omega(x,t)=\int_M H(x,y,
t)\Omega(y)\, dv_y
>0$ for $t>0$. By (2.5) we know that
$$
\int_1^t\Omega(x,s)\, ds \ge C\log t \tag 2.11
$$
for $t>>1$, for some $C>0$ independent of $t$. On the other hand,
by Theorem 3.1 of \cite{N1}, (2.9) implies that
$$
\int_1^t\Omega (x,s)\, ds \le \e \log t \tag 2.12
$$
for $t>>1$, in the {\it Case 1}. This proves the $(L, h)$ is flat
in the first case.

For the second case since $\tv_t(x,t)=\Omega(x,t)$, (2.11) implies
that
$$
\tv(x,t)\ge C\log t +C', \tag 2.13
$$
for $t>>1$ with positive constants $C$ and $C'$ independent of
$t$.
 By the assumption (2.4) we know that for any $\e>0$
$$
\aint_{B_o(r)}u\, dv \le \e\log r
$$
for $r>>1$. Using Theorem 3.1 of \cite{N1} we have that
$$
\tv(x,t)\le C(n)\e \log t
$$
for $t>>1$. This is a contradiction to (2.13).

The contradictions show that $(L, h)$ is flat in both cases. For
the last part of the theorem,  since $(L, h)$ is flat, $u$ is
harmonic. Then $u$ is a constant by a gradient estimate  of Cheng
and Yau \cite{CY}. We include here another proof based on the mean
value ineqality of Li-Schoen, which can be proved using the Moser
iteration argument even in the case gradient estimate fails. In
fact, for any $\a>1$, $\D u^\a =\frac{4(\a-1)}{\a}|\nabla
u^{\frac{\a}{2}}|^2$. Multiplying a cut-off function $\phi^2$ with
support in $B_o(2r)$ and integrating by parts we have
$$
\split
 2\frac{\a-1}{\a}\int_M|\nabla u^{\frac{\a}{2}}|^2\phi^2\, dv &=
 -\a\int_M u^{\a-1}<\nabla u, \nabla \phi>\phi\, dv\\
 &\le 2\int_M \lf(|\nabla u^{\frac{\a}{2}}|\phi\ri)\lf(
u^{\frac{\a}{2}}|\nabla \phi|
 \ri)\, dv\\
 &\le \frac{\a-1}{\a}\int_M|\nabla
 u^{\frac{\a}{2}}|^2\phi^2\, dv+\frac{\a}{\a-1}\int_M |\nabla \phi|^2
 u^{\a}\, dv.
\endsplit
$$
Therefore,
$$
\int_{B_o(r)}|\nabla u^{\frac{\a}{2}}|^2 \, dv\le
\frac{C_1}{r^2}\lf(\frac{\a}{\a-1}\ri)^2 \int_{B_o(2r)}u^{\a}\, dv
$$
for some universal constant $C_1$. Thus
$$
\aint_{B_o(r)}|\nabla u^{\frac{\a}{2}}|^2\le
C_2\lf(\frac{\a}{\a-1}\ri)^2\frac{\lf(\log r\ri)^{\a}}{r^2}\to 0,
\text { as } r\to \infty.
$$
Here $C_2$ is the constant independent of $r$.  On the other hand
since $|\nabla u|$ is subharmonic, the mean value inequality of
Li-Schoen \cite{LS} implies that for some $C_3(n)>0$
$$
\split \sup_{B_o(\frac{r}{2})}|\nabla u|^2& \le
C_3\aint_{B_o(r)}|\nabla u|^2\, dv\\
&\le C_2C_3\frac{(\log r)^{2}}{r^2}.
\endsplit
$$
Taking $r\to \infty$ we have that $|\nabla u|=0$. Thus  $u$ is a
constant.

\enddemo

\proclaim{Remark} The results in Theorme 2.2 were proved earlier
by Luen-Fain Tam and the author in \cite{NT2}. (cf. Theorem 0.3
and Corollary 6.1 therein.) The proof in \cite{NT2, Theorem 0.3}
uses the $L^2$-estimate of $\dbar$-operator of H\"ormander, as
well as Hamilton's strong maximum principle for tensors satisfying
(1.1). The proof of the gap result in \cite{NT2, Crollary 6.1}
uses the Liouville result above, along with the previous
developed,  quite sophisticated techniques of solving the
Poincar\'e-Lelong equation in \cite{MSY} and \cite{NST} as well as
some new refinements through heat equation deformation (cf.
Section 6 of \cite{NT2}).
\endproclaim

\proclaim{Corollary 2.2} Let $M$ be a complete K\"ahler manifold
with nonnegative holomorphic bisectional curvature. Let $u(x)$ be
a continuous plurisubharmonic functions of exponential growth. Let
$\tv(x,t)$ be the solution to $\heat \tv(x,t)=0$. Then
$$
w_t+\nabla_{\a}w\, V_{\bar{\a}}+\nabla_{\bar{\a}}w\, V_\a
+v_{\abb}V_{\bar{\a}}V_{\beta} +\frac{w}{t}\ge 0. \tag 2.14
$$
holds on $M\times (0, \infty)$ for any $(1,0)$ vector field $V$.
Here $w(x,t)=\Delta \tv (x,t)$.
\endproclaim
\demo{Proof} It follows from Theorem 1.1, Theorem 2.1 and Lemma
2.1.
\enddemo

\proclaim{Remark} The estimate of the form (2.14) was first proved
for the plurisubharmonic functions deformed by the time-dependent
(with metric evolved by K\"ahler-Ricci flow) heat equation in
\cite{NT1}. It was also used to prove the Liouville theorem for
the plurisubharmonic fcuntions for the first time there. However,
due to complications caused by the K\"ahler-Ricci flow, the result
requires various assumptions on the curvature of the initial
metric on $M$. In particular, one has to assume the boundedness of
the curvature, which is rather artificial for the study of the
function theory on $M$.
\endproclaim

\input amstex
\documentstyle{amsppt}
\magnification=1200 \hsize=13.8cm \catcode`\@=11
\def\NoLogo{\let\logo@\empty}
\catcode`\@=\active \NoLogo
\def\heatt{\lf (\Delta-\frac{\p}{\p t}\ri)}

\def\heat{\lf(\frac{\p}{\p t}-\Delta\ri)}

\def \b {\beta}

\def\lf{\left}
\def\ri{\right}
\def\bbar{\bar \beta}
\def\a{\alpha}

\def\g{\gamma}
\def\e{\epsilon}
\def\p{\partial}
\def\delbar{{\bar\delta}}
\def\ddbar{\partial\bar\partial}
\def\dbar{\bar\partial}

\def\C{\Bbb C}
\def\R{\Bbb R}
\def\N{\Bbb N}

\def\vp{\varphi}

\def\dbar{\bar\partial}

\def\abb{{\alpha\bar\beta}}
\def\gbd{{\gamma\bar\delta}}

\def \D {\Delta}
\def\aint{\frac{\ \ }{\ \ }{\hskip -0.4cm}\int}
\vsize=19.0 cm

\subheading{\S3 Dimension estimates I}

\vskip .2cm

Let $M$ be a complete K\"ahler manifold with nonnegative
bisectional curvature of complex dimension $m$. We shall show
further applications of gradient estimate (2.3) in this section in
the study of  holomorphic functions of polynomial growth. Let us
first fix the notation. We call a holomorphic function $f$ of
polynomial growth if there exists $d\ge 0$ and $C=C(d, f)$ such
that
$$
|f|(x)\le C(r^d(x)+1), \tag 3.1
$$
where $r(x)$ is the distance function to a fixed point $o\in M$.
For any $d>0$ we denote ${\Cal O}_d(M)=\{f\in {\Cal O}(M)| f(x)
\text{ satisfies } (3.1)\}.$
 Let ${\Cal O}_P(M)$ denotes the space of the
holomorphic functions of polynomial growth. Since any sub-linear
growth holomorphic function is constant on a complete K\"ahler
manifold with nonnegative Ricci curvature,
$$
{\Cal O}_P(M)=\C \cup \lf(\cup_{d\ge 1} {\Cal O}_d(M)\ri).
$$
We also define the order of $f$ in the sense of Hadamard to be
$$
\text{Ord}_{H}(f) =\limsup_{r\to \infty} \frac{\log \log
A(r)}{\log r}
$$
where $A(r)=\sup_{B_o(r)}|f(x)|$. It is clear that if $f\in {\Cal
O}_P(M)$, $Ord_H(f)=0$. We call $f$ has finite order if
$Ord_H(f)<\infty$.  The first issue we are going to address is
estimating the dimension of ${\Cal O}_d(M)$. Let us start with
some simple observations.

\proclaim{Lemma 3.1} Let $f\in {\Cal O}(M)$ be a non-constant
holomorphic functions of finite order. Denote $u(x)=\log
(|f|(x))$. Then there exist solution $v(x,t)$ to $\heatt v(x,t)=0$
such that $v(x,0)=u(x)$, $v(x,t)$ is plurisubharmonic. Moreover,
the function  $w(x,t)=\D v(x,t)>0$, for $t>0$, and
$$
\frac{\p}{\p t}\lf(t\, w(x,t)\ri)\ge 0 \tag 3.2
$$
\endproclaim
\demo{Proof} Let $u_j(x)=\log (|f|(x)+\frac{1}{j})$. Then
$$
v_j(x,t)=\int_M H(x,y,t)u_j(y)\, dy \tag 3.3
$$
gives the solution $v_j(x,t)$ such that $v_j(x,0)=u_j(x)$. Clearly
$v_j(x,t)$ satisfies assumption of Lemma 2.1. Thus $v_i(x,t)$ are
plurisubharmonic functions.  Let $j\to \infty$ in (3.3), we obtain
$$
v(x,t)=\int_M H(x,y, t)u(y)\, dy
$$
a solution with $v(x,0)=u(x)$.  Let $w(x,t)=v_t(x,t)$. Since
$\{v_j\}$ is a decreasing sequence $v(x,t)$ is also
plurisubharmonic. To prove (3.2) we claim that
$w_j(x,t)=(v_j)_t(x,t)$ satisfies (3.2). Since $w_j(x,t)\to
w(x,t)$ uniformly on compact subsets of $M\times(0, \infty)$ the
claim implies that $w(x,t)$ also satisfies (3.2). In order to
prove (3.2) for $w_j$, we notice that $w_j(x,0)=\D u_j(x)$. By the
strong maximum principle we have that $w_j(x,t)>0$, otherwise
$u_j$ is harmonic, which implies that $f$ is a constant by
Cheng-Yau's gradient estimate \cite{CY}. This proves (3.2). To
show that $w(x,t)>0$, observe that $\lim_{t\to 0}w(x,t)=\D u(x)$.
We claim that $f$ must vanish at some place. Otherwise $u$ is a
harmonic function of sub-linear growth, which implies $u$ is a
constant by Cheng-Yau's gradient estimate \cite{C-Y} again. This
then implies $f$ is a constant, which contradicts with the
assumption. Therefore $\D u$ must be a non-zero, nonnegative
measure.  This implies that $w(x,t)$ can not be identically zero
for $t>0$. By the strong maximum principle we then have that
$w(x,t)>0$ for $t>0$.
\enddemo

\proclaim{Remark} On a complete K\"ahler manifold $M$ with
nonnegative Ricci curvature, since Cheng-Yau's gradient estimate
implies that any sub-linear growth harmonic function is constant,
it then implies that any non-constant holomorphic function $f$
with $Ord_H(f)<1$ must vanishes somewhere. This in particular
generalizes the fundamental theorem of algebra to the complete
K\"ahler manifolds with nonnegative Ricci curvature.
\endproclaim

Recall that for any positive $(p, p)$ current $\Theta$  one can
define the Lelong number of  $\Theta$ at $x$ as
$$
\nu(\Theta, x)=\lim_{r\to 0}\nu(\Theta, x, r) \tag 3.4
$$
where
$$
\nu(\Theta, x,
r)=\frac{1}{r^{2(m-p)}\pi^{m-p}}\int_{B_x(r)}\Theta\wedge
\omega^{m-p}. \tag 3.5
$$
The existence of the limit in (3.4) is ensured by (0.4).
 For $f(x)\in {\Cal O}(M)$ we denote $Z_f$ to be the zero set
of $f$. $Z_f$ is a positive $(1,1)$ current. The Poincar\'e-Lelong
Lemma states that
$$
\frac{\sqrt{-1}}{2\pi}\ddbar \log (|f|^2)=Z_f. \tag 3.6
$$
We define
$$ord_x(f)=\max\{ m\in \N | D^{\a} f(x)=0, |\a|<m\}.$$
It is well-know that
$$
ord_x(f)=\nu(Z_f, x). \tag 3.7
$$
 One can refer \cite{D} or \cite{GH} for
details of above cited results on the Lelong number and
$ord_x(f)$. Using (3.4)--(3.7) one can have that
$$
ord_x(f)=\frac{1}{2m}\lim_{r\to
0}\lf(\frac{r^2}{V_x(r)}\int_{B_x(r)}\D \log |f|\, dv\ri). \tag
3.8
$$

\proclaim{Theorem 3.1} Let $M$ be a complete K\"ahler manifold
with nonnegative holomorphic bisectional curvature of complex
dimension $m$. Then there exists a constant $C_1=C_1(m)$ such that
for any $f\in {\Cal O}_d(M)$
$$
ord_x(f)\le C_1 d. \tag 3.9
$$
In particular, it implies
$$
\dim_{\C}({\Cal O}_d(M))\le C_2 d^m \tag 3.10
$$
for some $C_2=C_2(m)$.
\endproclaim
\demo{Proof} Let $u(x), \  v(x,t), $  and $w(x,t)$ be as in Lemma
3.1. From (3.2) we know that
$$
\lf(tw(x,t)\ri)_t\ge 0. \tag 3.11
$$
We are going to show that there exist positive constants
$C_3=C_3(m)$ and $C_4=C_4(m)$ such that
$$
\lim_{t\to 0}(tw(x,t))\ge C_3\,  ord_x(f) \tag 3.12
$$
and
$$
tw(x,t)\le  C_4\, d \tag 3.13
$$
for $t>>1$.

We first show (3.12). The approximation argument in the proof of
Lemma 3.1 shows that
$$
w(x,t)=\int_M H(x,y, t)\D \log(|f|(y))\, dv_y.
$$
As in the proof of Theorem 3.1 in \cite{N1}, using the Li-Yau's
\cite{LY} lower bound of heat kernel we have that
$$
w(x,t)\ge C(m) \frac{1}{V_x(\sqrt{t})}\int_{B_x(\sqrt{t})}\D \log
(|f|(y))\, dv_y.
$$
Therefore
$$
tw(x,t)\ge C(m)\frac{t}{V_x(\sqrt{t})}\int_{B_x(\sqrt{t})}\D \log
(|f|(y))\, dv_y.
$$
Now (3.12) follows easily from (3.8).

To prove (3.13), we first observe that,  by Theorem 3.1 in
\cite{N1} , for $t>>1$
$$
v(x,t)\le C_5 d \log t \tag 3.14
$$
for some constant $C_5=C_5(m)$,  since from the assumption (3.1)
one has $\log |f|(x) \le d\log (r(x)+1) +C$. (Here one can not
apply Theorem 3.1 of \cite{N1} directly since $v$ is not always
nonnegative. But we can use $u_{+}$ as the initial date to obtain
a solution to the heat equation, which serves a barrier from above
for $v$.) We claim that this implies
$$
tw(x,t)\le C_5 \, d
$$
for $t>>1$.  Otherwise, we have some $\e>0$ such that
$$
tw(x,t)>(C_5+\e) d
$$
for $t>>1$. Here we have used the monotonicity of $tw(x,t)$.
Therefore
$$
v(x,t)\ge (C_5+\e)d\log t -A
$$
where $A$ is independent of $t$. This contradicts  (3.14). Since
(3.9) follows from (3.11)--(3.13) and (3.10) follows from (3.9) by
a simple dimension counting argument (cf. \cite{M, page 221}) we
complete the proof of the theorem.
\enddemo

\proclaim{Remark} The dimension estimate as well as the
multiplicity estimate (3.9) for the holomorphic functions of
polynomial growth was first considered in \cite{M} by Mok. In
\cite{M}, the estimate was obtained for manifolds with maximum
volume growth as well as  a point-wise quadratic decay assumption
on the curvature (Cf. (0.5) and (0.6)). Also the constant in the
estimate similar to (3.9), obtained in \cite{M}, depends on the
local geometry of $M$. Here the constant depends only on the
complex dimension. The estimate (3.10) is sharp in the power.
\endproclaim

Denote ${\Cal M}(M)$  the function field generated by ${\Cal
O}_P(M)$. Namely any $F\in {\Cal M}(M)$ can be written as
$F=\frac{g}{h}$ with $ g, h\in {\Cal O}_P(M)$. A direct
consequence of Theorem 3.1 is the following statement.

\proclaim{Corollary 3.1} Let $M$ be as in Theorem 3.1. Then the
transcendence  degree of ${\Cal M}(M)$ over $\C$, $deg_{tr}({\Cal
M}(M))$ satisfies
$$
deg_{tr}({\Cal M}(M))\le m.
$$
Moreover, in  the equality case, ${\Cal M}(M)$ is a finite
algebraic  extension over $\C(f_1, \cdots,f_m)$, where $f_i$ are
the transcendental elements in ${\Cal M}(M)$. More precisely,
there exist $g, h \in {\Cal O}_P(M)$ and  a polynomial $P$ with
coefficients in $\C(f_1, \cdots,f_m)$  such that
$P(\frac{g}{h})=0$ and ${\Cal M}(M)=\C(f_1, \cdots, f_m,
\frac{g}{h}).$
\endproclaim
\demo{Proof} It follows from the so-called Poincar\'e-Siegel
arguments. See, for example \cite{M, pages 220-221} or \cite{S,
pages 176-178}.
\enddemo

The following result can be viewed as a gap theorem for
holomorphic functions. It is related to Theorem 0.3 of \cite{NT2}.

\proclaim{Corollary 3.2} Let $M$ be as in Theorem 3.1. There
exists a $\e=\e(m)>0$ such that
$$
\dim({\Cal O}_{1+\e}(M))\le m+1.
$$
\endproclaim
\demo{Proof} By Lemma 4.4 of \cite{NT2}, we know that if $f\in
{\Cal O}(M)$, $\log |\nabla f|$ is plurisubharmonic. Also, by
Cheng-Yau's gradient estimate, $|\nabla f|\le C(r^{\e}+1)$. Apply
the proof of Theorem 3.1 we know that
$$
ord_x(Df) \le C_1(m)\e.
$$
Here $Df=(\frac{\p f}{\p z_1}, \frac{\p f}{\p z_2}, \cdots,
\frac{\p f}{\p z_m})$ with respect to a fixed local coordinates
chart. If we choose $\e=\frac{1}{2C_1}$, we have that
$$
ord_x(Df)\le \frac{1}{2}.
$$
This implies that $Df(x)\ne 0$. The dimension counting argument
then gives the conclusion.
\enddemo

The example of `round-off' cones on pages 3--4 of \cite{NT2} shows
that one can not expect that ${\Cal O}_{1+\e}(M)={\Cal O}_1(M)$.
Namely, one can not conclude that $f$ is indeed linear. However,
this is the case if one assume stronger `closeness' assumption as
in Theorem 0.3 of \cite{NT2}.

The  dimension estimates for the holomorphic functions can be
generalized for the holomorphic sections of polynomial growth of
holomorphic line bundles with controlled positive part of the
curvature. In particular, it applies to the non-positive line
bundles. (We call $(L, h)$ is non-positive if $\Omega_{\abb}(x)\le
0$.) We treat the non-positive line bundle with continuous
curvature in this section first and leave the more complicated
case when the curvature has positivity to the next section. Before
we state the result let us denote
$$
{\Cal O}_{d}(M, L)=\{s\in {\Cal O}(M, L)\, |\,  \|s\|(x)\le
C(r(x)+1)^d\}.
$$
Here $r(x)$ is the distance function to a fixed point $o\in M$.

\proclaim{Theorem 3.2} Let $M$ be a complete K\"ahler manifolds
with nonnegative bisectional curvature.  Let $(L, h)$ be a
hermitian line bundle with non-positive curvature.  Then
$$
\dim({\Cal O}_{d}(M, L))\le C_1 d^m. \tag 3.15
$$
Here $C_1=C_1(m)$.
\endproclaim
\demo{Proof} We assume that there exists $s\in {\Cal O}_{d}(M,
L)$. The well-known Poincar\'e-Lelong equation states
$$
\frac{\sqrt{-1}}{2\pi}\ddbar \log (\|s\|^2)=[s]-\rho. \tag 3.16
$$
Here $\rho=\frac{\sqrt-1}{2\pi}\Omega_{\abb}dz_{\a}\wedge
d\bar{z}_{\beta}$, $[s]$ is the divisor defined by the zero locus
of $s$. In particular, it implies that
$$
\D \log \|s\|^2(x) \ge -\Omega(x). \tag 3.17
$$
Now let $u(x)=\log(\|s\|)$ and solve the heat equation $\heatt
v(x,t)=0$ with the initial data $v(x,0)=u(x).$ The solvability can
be justifies by the argument of Lemma 3.1. Similarly we  have that
$v(x,t)$ is plurisubharmonic and $w(x,t)=v_t(x,t)$ satisfying
$\heatt w(x,t)=0$, $w(x,0)=\D \log \|s\|$. Moreover,
$$
\lf(tw(x,t)\ri)_t \ge 0.
$$
The argument of Lemma 3.1 also implies
$$
w(x,t)=\int_M H(x,y, t)\lf(\D \log \|s\|(y)\ri) dv_y. \tag 3.18
$$
We denote $mult_{x}([s])$ the multiplicity of the divisor $[s]$.
It is from the definition that
$$
mult_x([s])=\frac{1}{2m}\lim_{r\to
0}\lf(\frac{r^2}{V_x(r)}\int_{B_x(r)}\D \log \|s\|(y)\, dv_y\ri).
$$
Now the same argument as the proof  of  (3.12) shows that
$$
\lim_{t\to 0}tw(x,t)\ge C_2(m)\, mult_x([s]). \tag 3.19
$$
We claim that
$$
\lim_{t\to \infty} t w(x,t) \le C_3(m) d. \tag 3.20
$$
In fact since $v(x,t)=\int_M H(x,y, t) \log \|s\|(y)\, dv_y$ we
have that
$$
v(x,t)\le C_4(m) d\log t \tag 3.21
$$
for $t>>1$,  by Theorem 3.1 of \cite{N1}, as in the proof of
Theorem 3.1. Now similar argument as in Theorem 3.1 shows (3.20).
Therefore we have
$$
mult_x([s])\le C_5(m) d, \tag 3.22
$$
from which (3.15) follows by dimension counting argument.
\enddemo
The proof of the above result as well as the proof of Theorem 3.1
gives the following improvement of a earlier result \cite{Theorem
4.3, NT2}.

\proclaim{Corollary 3.3} Let $M$ be a complete K\"ahler manifold
with nonnegative bisectional curvature. Assume that $M$ admits a
nonconstant holomorphic function of polynomial growth and the
bisectional curvature is positive at some point. Then
$$
V_x(r)\ge \frac{C_3}{r^{m+1}} \tag 3.23
$$
and
$$
\aint_{B_x(r)}{\Cal R}(y)\, dv_y \le \frac{C_4}{r^2}\tag 3.24
$$
for some positive constant $C_3$ and $C_4$ (which might depends on
$x$).
\endproclaim
\demo{Proof} We only prove (3.24) here and leave (3.23) to the
interested reader. By the assumption that $M$ admits a holomorphic
functions of polynomial growth and $M$ has quasi-positive
bisectional curvature, the proof of \cite{Theorem 4.3, NT2}
implies that there exists a smooth strictly plurisubharmonic
function $u(x)$ on $M$ such that $u(x)\le C(\log r(x)+2)$. By the
proof of Lemma 4.2 of \cite{NT2} we can have a nontrivial $s\in
{\Cal O}_d(M, K_M)$ for some $d>0$. Now we apply the argument of
Theorem 3.2 to the case $L=K_M$ and have that $-\Omega(y)={\Cal
R}(y)$. Now combining (3.17) and (3.18), we have that
$$
w(x,t)\ge \int_M H(x,y, t){\Cal R}(y)\, dv_y.
$$
Applying Theorem 3.1 of \cite{N1} we then have
$$
w(x,t)\ge C_5(m) \aint_{B_x(\sqrt{t})}{\Cal R}(y)\, dv_y. \tag
3.25
$$
Now (3.24) follows from (3.20), (3.25) and the monotonicity of
$tw(x,t)$.
\enddemo

Note that the special case $m=1$ of (3.23) recovers the earlier
result of Wu \cite{W}.

\proclaim{Remark} Theorem 3.2 can be applied to the canonical line
bundle to give the dimension estimates for the canonical sections
of polynomial growth. In Corollary 3.3, (3.24) also holds for the
case with general non-positive line bundles by assuming that there
exists a holomorphic section of polynomial growth. In general, the
estimate holds for the negative part of $\Omega$.

We conjecture that under the assumption of Corollary 3.3 one
should be able to prove  that $M$ has maximum volume growth. The
intuition for this is that it seems that every transcendental
 holomorphic function of polynomial growth seems to contribute to
 the
 volume by a factor of $r^2$. On the other hand, under the assumption of
 Corollary 3.3 one in fact has  $(f_1, \cdots, f_m)$ to form a local
 coordinate near  any given point.
\endproclaim

\input amstex
\documentstyle{amsppt}
\magnification=1200 \hsize=13.8cm \catcode`\@=11
\def\NoLogo{\let\logo@\empty}
\catcode`\@=\active \NoLogo
\def\heatt{\lf (\Delta-\frac{\p}{\p t}\ri)}

\def\heat{\lf(\frac{\p}{\p t}-\Delta\ri)}

\def \b {\beta}

\def\lf{\left}
\def\ri{\right}
\def\bbar{\bar \beta}
\def\a{\alpha}

\def\g{\gamma}
\def\e{\epsilon}
\def\p{\partial}
\def\delbar{{\bar\delta}}
\def\ddbar{\partial\bar\partial}
\def\dbar{\bar\partial}

\def\C{\Bbb C}
\def\R{\Bbb R}
\def\N{\Bbb N}
\def\tv{\tilde v}
\def\vp{\varphi}

\def\dbar{\bar\partial}

\def\abb{{\alpha\bar\beta}}
\def\gbd{{\gamma\bar\delta}}

\def \D {\Delta}
\def\aint{\frac{\ \ }{\ \ }{\hskip -0.4cm}\int}
\vsize=19.0 cm

\subheading{\S4 Dimension estimates--the sharp ones}

\vskip .2cm

Let $M$ be a complete K\"ahler manifold of complex dimension $m$.
Under the assumption $M$ has nonnegative Ricci curvature, the
function $\frac{V_x(r)}{r^n}$ is monotone decreasing ($n=2m$ is
the real dimension). If it has a positive limit $\theta=\lim_{r\to
\infty}\frac{V_x(r)}{r^n}$ we call $M$ is of maximum volume
growth.
 In \cite{LW} the authors proved some  asymptotically
sharp dimension estimates for the harmonic functions of polynomial
growth on a complete Riemannian manifold with nonnegative
sectional curvature and maximum volume growth.  Here we shall show
 the sharp dimension estimate for ${\Cal
O}_d(M)$ for complete K\"ahler manifolds with nonnegative
bisectional curvature and maximum volume growth.

\proclaim{Theorem 4.1} Let $M^m$ be a complete K\"ahler manifold
with nonnegative holomorphic bisectional curvature. Assume that
$M$ is of maximum volume growth. Then
$$
ord_x(f)\le [d]. \tag 4.1
$$
In particular,
$$
\dim_{\C}({\Cal O}_d(M))\le \dim_{\C}({\Cal O}_{[d]}(\C^m)). \tag
4.2
$$
Here $[d]$ is the greatest integer less than or equal to $d$.
\endproclaim

We need several lemmas to prove the above results. The first one
is the sharpen version of (3.12).

\proclaim{Lemma 4.1} Let $u, v, w$ be as in Lemma 3.1. Then
$$
\lim_{t\to 0}tw(x,t)=\frac{1}{2}ord_x(f). \tag 4.3
$$
\endproclaim
\demo{Proof} Since $w(x,t)$ solves the heat equation $\heatt
w(x,t)=0$ with the initial data being the positive measure $\D
\log (|f|(y))$, we can apply Theorem 3.1 of \cite{N1} to $w$. (It
is easy to see that the result proved in \cite{N1} can be
generalized to the case when the initial dada being the positive
measure.) Therefore we have that
$$
\aint_{B_x(r)}\D \log(|f|(y))\, dv_y \le C(m) w(x, r^2) \tag 4.4
$$
for some constant $C(m)>0$. On the other hand, from the proof of
Theorem 3.1 we know that
$$
tw(x,t)\le C_5\, d \tag 4.5
$$
for $t>>1$,  which then implies
$$
\aint_{B_x(r)}\D \log(|f|(y))\, dv_y \le  \frac{C(m)\,d}{r^2} \tag
4.6
$$
for $r>>1$.

It is well know that
$$
H(x,y, t)\sim \frac{1}{(4\pi
t)^{\frac{n}{2}}}\exp(-\frac{r^2(x,y)}{4t}) +\text{ lower order
term} \tag 4.7
$$
as $t\to 0$. By (3.8) we also know that  for $\e>0$ there exits
$\delta>0$ such that
$$
2m \, ord_x(f)-\e \le \frac{r^2}{V_x(r)}\int_{B_x(r)}\D \log |f|\,
dv\le 2m\,  ord_x(f)+\e \tag 4.8
$$
for $r\le \delta$. Write
$$
\split tw(x,t) &=t\int_{M} H(x,y,t)\D \log(|f|(y))\, dv_y\\
&=t\int_{r(x,y)\ge \delta}H(x,y,t)\D \log(|f|(y))\,
dv_y+t\int_{r(x,y)\le \delta} H(x,y,t)\D \log(|f|(y))\, dv_y\\
&= I +II.
\endsplit
$$
Here $I$ and $II$ denote the first and the second term in the
second line respectively. In the following we are going to show
that $I$ has limit $0$, as well as
$$
\frac{1}{2}\, ord_x(f)-2\e \le \liminf_{t \to 0}II \tag 4.9
$$
and
$$ \limsup_{t\to 0}II \le \frac{1}{2}\, ord_x(f)+2\e. \tag 4.10
$$
Clearly, (4.1) is a consequence of these conclusions. Using
Li-Yau's upper bound on the heat kernel estimate we have
$$
\split I &\le \frac{C(n)t}{V_x(\sqrt{t})}\int_\delta
^{\infty}\exp(-\frac{s^2}{5t})\lf(\int_{\p B_x(s)}\D \log
(|f|(y))\, dA\ri)\, ds\\
&\le
\frac{C(n)t}{V_x(\sqrt{t})}\exp(-\frac{\delta^2}{5t})\int_{B_x(\delta)}\D
\log |f|(y)\, dv_y \\
&\ \  + C(n)t\int_{\delta}^{\infty}
\exp(-\frac{s^2}{5t})\lf(\frac{s}{\sqrt{t}}\ri)^{n}\lf(\aint_{B_x(s)}\D
\log (|f|(y))\, dv_y\ri)\lf(\frac{2s}{5t}\ri)\, ds\\
&= III+ IV.
\endsplit
$$
Here we have used the volume comparison theorem and assumed that
$\sqrt{t}\le \delta$. $III$ and $IV$ denote the term in the second
and third line respectively. Clearly $\lim_{t\to 0} III=0$. On the
other hand
$$
IV \le C(n)\int_{\frac{\delta^2}{5t}}^{\infty}\exp(-\tau)
\tau^{\frac{n}{2}-1} t\tau \lf(\aint_{B_x(\sqrt{5t\tau})}\D \log
|f|(y)\, dv_y\ri)\, d\tau
$$
Using the estimate (4.6) we have that $\lim_{t\to 0} IV =0$.
Therefore we have shown that
$$
\lim_{t\to 0} I =0. \tag 4.11
$$
Now we prove (4.10).  Using (4.7),  for $t<<1$,
$$
\split II &\le t\int_0^{\delta}\frac{1}{(4\pi
t)^{\frac{n}{2}}}\exp(-\frac{s^2}{4t})\lf(\int_{\p B_x(s)}\D
\log|f|(y)\, dy\ri)\, ds +\frac{\e}{2}\\
& = t \frac{1}{(4\pi
t)^{\frac{n}{2}}}\exp(-\frac{\delta^2}{4t})\lf(\int_{
B_x(\delta)}\D \log|f|(y)\,
dy\ri)\\
&\ \ +t\int_0^{\delta}\frac{1}{(4\pi
t)^{\frac{n}{2}}}\exp(-\frac{s^2}{4t})\lf(\int_{B_x(s)}\D
\log|f|(y)\, dy\ri)\lf(\frac{s}{2t}\ri)\, ds +\frac{\e}{2}\\
&=V +VI +\frac{\e}{2}.
\endsplit
$$
Here $V$ and $VI$ denote the term in the second and the third line
respectively. The  term $V$ has limit $0$ as in the estimate of
$I$. To estimate the  term $VI$ we use (4.8) and the fact that
$$
V_x(s) \sim \omega_{n} s^n
$$
for $s\to 0$, where $\omega_n$ is the volume of unit ball in
$\R^n$. In deed,
$$
\split VI &\le t\int_0^{\delta}\frac{\omega_n s^n}{(4\pi
t)^{\frac{n}{2}}}\exp(-\frac{s^2}{4t})\lf(\aint_{B_x(s)}\D
\log|f|(y)\, dy\ri)\lf(\frac{s}{2t}\ri)\, ds\\
&=\frac{\omega_n}{\pi^{\frac{n}{2}}}\int_0^{\frac{\delta^2}{4t}}\frac{1}{4}\exp(-\tau)
\tau^{\frac{n}{2}-1}\lf( (4t\tau)\aint_{B_x(\sqrt{4t\tau})}\D \log
|f|(y)\, dv_y\ri)d \tau \\
&\le \frac{m}{2}\frac{\omega_n}{\pi^{\frac{n}{2}}} \,
ord_x(f)\int_0^{\frac{\delta^2}{4t}}\exp(-\tau)
\tau^{\frac{n}{2}-1}d\tau +\frac{\e}{2}\\
&=\frac{m}{2}\frac{\omega_n}{\pi^{\frac{n}{2}}} \, ord_x(f)
\Gamma(\frac{n}{2})+\frac{\e}{2}\\
&=\frac{1}{2}\, ord_x(f)+\frac{\e}{2}.
\endsplit
$$
This proves (4.10). The proof for (4.9) is similar.
\enddemo

Notice that we do not make use of the maximum volume growth in
Lemma 4.1. The next lemma sharpens (3.13), which makes use of the
maximum volume growth assumption.

\proclaim{Lemma 4.2} Let $M$ be as in Theorem 4.1 and let $u, v,
w$ be as in Lemma 3.1. Then
$$
\limsup_{t\to \infty} \frac{v(x,t)}{\log t}\le \frac{1}{2}d. \tag
4.12
$$
which, in particular,  implies
$$ \lim_{t\to \infty} t w(x,t)\le
\frac{1}{2} d. \tag 4.13
$$
\endproclaim
In order to prove the above lemma we need a result of Li-Tam-Wang
\cite{LTW, Theorem 2.1} on the upper bound of the heat kernel
under the maximum volume growth assumption.

\proclaim{Theorem 4.2 (Li-Tam-Wang)} Let $M$ be a complete
Riemannian manifold with nonnegative Ricci curvature and maximum
volume growth. For any $\delta>0$, the heat kernel of $M$ must
satisfy the estimate
$$
\split \frac{\omega_n}{\theta   (\delta r(x,y))}&\, (4\pi
t)^{-\frac{n}{2}}\exp\lf( -\frac{1+9\delta}{4t}r^2(x,y)\ri)\le
H(x,y,t)\\
&\le \lf(1+C(n,\theta)(\delta
+\beta)\ri)\frac{\omega_n}{\theta}(4\pi
t)^{\frac{n}{2}}\exp\lf(-\frac{1-\delta}{4t}r^2(x,y)\ri),
\endsplit \tag 4.14
$$
where
$$
\beta =\delta^{-2n}\max_{r\ge
(1-\delta)r(x,y)}\{1-\frac{\theta_x(r)}{\theta_x(\delta^{2n+1}r)}\}.
\tag 4.15
$$
\endproclaim
Note that $\beta$ is a function of $r(x,y)$ such that
$$
\lim_{r(x,y)\to \infty} \beta =0. \tag 4.16
$$

\demo{Proof of Lemma 4.2} By (4.16) and the fact that
$$
A_x(s) \sim n s^{n-1}\theta,
$$
as $s\to \infty$, where $A_x(s)$ is the area of the sphere $\p
B_x(s)$, we know that for any $\e>0$, there exists a positive
constant $A>0$ such that for $s=r(x,y)\ge A$
$$
\beta \le \e \text{  and } \frac{A_x(s)}{\theta}\le \lf((1+\e)
n\ri)s^{n-1}. \tag 4.17
$$
Now, we estimate $v(x,t)=\int_M H(x,y,t)\log|f|(y)\, dv_y$. Using
the upper bound of Li-Yau we have that
$$
\split
 v(x,t)  & = \int_{r(x,y)\le A} H(x,y, t) \log |f|(y)\, dv_y +\int_{r(x,y)\ge A} H(x,y, t)
 \log |f|(y)\, dv_y \\
&\le \frac{C(n)}{V_x(\sqrt{t})}\int_0^A
\exp(-\frac{s^2}{5t})\lf(\int_{\p B_x(s)}\log |f|(y)\, dA_y\ri)\,
ds +II\\
&\le \frac{C(n)}{V_x(\sqrt{t})}\int_{B_x(A)}\log |f|(y)\, dv_y
+II\\
&= I +II.
\endsplit
$$
Here we use  $II$ to represent the second term of the first line
and $I$ to represent the first term of the third line. Clearly
$$
\lim_{t\to \infty} \frac{I}{\log t}=0. \tag 4.18
$$
We claim that
$$
\limsup_{t\to \infty} \frac{II}{\log t} \le \frac{1}{2}\, d .\tag
4.19
$$
The lemma follows easily from (4.18) and (4.19).  To prove (4.19)
we need the estimate (4.14) of Li-Tam-Wang. Notice that this is
the only place the maximum volume growth condition is used. In
deed, by (4.14) and (4.17), for the given fixed $\delta $ and
$\e>0$,
$$
\split II&\le\lf( 1+C(n,\theta)(\delta
+\e)\ri)\frac{\omega_n}{(4\pi
t)^{\frac{n}{2}}}\int_A^{\infty}\exp(-\frac{1-\delta}{4t}s^2)\frac{1}{\theta}\lf(\int_{\p
B_x(s)}\log |f|(y)\, dA_y\ri)\, ds\\
&\le \lf( 1+C(n,\theta)(\delta +\e)\ri)\frac{\omega_n}{(4\pi
t)^{\frac{n}{2}}}\int_A^{\infty}\exp(-\frac{1-\delta}{4t}s^2)\frac{A_x(s)}{\theta}
\lf(d\log s +\tilde C\ri)\, ds\\
& \le \lf( 1+C(n,\theta)(\delta +\e)\ri)\frac{\omega_n}{(4\pi
t)^{\frac{n}{2}}}\int_A^{\infty}\exp(-\frac{1-\delta}{4t}s^2)n(1+\e)s^{n-1}
\lf(d\log s +\tilde C\ri)\, ds.
\endsplit
$$
Here $\tilde C$ is the constant in (3.1).  Let
$$
III=\lf( 1+C(n,\theta)(\delta +\e)\ri)(1+\e)\frac{d\, n \,
\omega_n}{(4\pi
t)^{\frac{n}{2}}}\int_A^{\infty}\exp(-\frac{1-\delta}{4t}s^2)s^{n-1}
\log s \, ds
$$
and
$$
IV=\lf( 1+C(n,\theta)(\delta +\e)\ri)\tilde C(1+\e)
\frac{n\omega_n}{(4\pi
t)^{\frac{n}{2}}}\int_A^{\infty}\exp(-\frac{1-\delta}{4t}s^2)s^{n-1}
\, ds.
$$
Then we have $II\le III+IV$.  It is easy to check that
$$
\lim_{t\to \infty} \frac{IV}{\log t}=0.
$$
For $III$ we have that
$$
III\le \lf( 1+C(n,\theta)(\delta +\e)\ri)\frac{(1+\e)d\, n \,
\omega_n}{4(\pi(1-\delta))^{\frac{n}{2}}}\int_{\frac{(1-\delta)A^2}{4t}}^{\infty}
\exp(-\tau)\tau^{\frac{n}{2}-1} \lf(\log t
+\log(\frac{4\tau}{1-\delta}\ri) \, d\tau.
$$
Then we have that
$$
\split \limsup_{t\to \infty} \frac{III}{\log t}& \le \lf(
1+C(n,\theta)(\delta +\e)\ri)\frac{(1+\e)d\, n \,
\omega_n}{4(\pi(1-\delta))^{\frac{n}{2}}}\int_{0}^{\infty}
\exp(-\tau)\tau^{\frac{n}{2}-1}\, d\tau\\
&=\lf( 1+C(n,\theta)(\delta +\e)\ri)\frac{(1+\e)d\, n \,
\omega_n}{4(\pi(1-\delta))^{\frac{n}{2}}}\Gamma(\frac{n}{2})\\
&=\lf( 1+C(n,\theta)(\delta
+\e)\ri)\frac{(1+\e)d}{2(1-\delta)^{\frac{n}{2}}}.
\endsplit
$$
Since $\delta $ and $\e$ are arbitrary chosen positive constants,
this proves (4.19).
\enddemo

\demo{Proof of Theorem 4.1} By Lemma 4.1 and 4.2 we have that
(4.1), from which the theorem follows by the dimension counting.
More precisely one can define the map
$$
\Phi: {\Cal O}_d(M) \to \C^{K_{[d]}}
$$
by for all $|\a|\le [d]$,
$$\Phi(f)=(f(x), Df(x), \cdots, D^\a f(x))
$$
where $K_{[d]}=\dim({\Cal O}_d(M))=1+\lf(\matrix
n\\n-1\endmatrix\ri)+\cdots +\lf(\matrix
n+[d]-1\\n-1\endmatrix\ri)=\lf(\matrix n+[d]\\ n\endmatrix\ri)$.
Assume that the conclusion of the theorem is not true. Then there
exists $f\ne 0$, such that $\Phi(f)=0$. This implies that
$ord_x(f)>d$, which is a contradiction to (4.1).
\enddemo

\proclaim{Remark} Theorem 4.1 was conjectured by Yau \cite{Y}. It
is unclear if the result holds for general case without assuming
the maximum volume growth. The only place where the maximum volume
growth assumption was used is in Lemma 4.2, where  the
asymptotically sharp heat kernel upper bound estimate of Li-Tam
Wang was applied. If the similar heat kernel estimate holds for
the general case without assuming the maximum volume growth the
argument  given above  gives the sharp dimension bound for general
case without much modification.

A localized version of the estimate  (4.1) was proved  for
polynomials on $\C^m$ by Bombieri \cite{B} in the study of the
algebraic values of meromorphic maps. Similar localization can be
derived from the above proof of Theorem 4.1 for the holomorphic
functions of polynomial growth on complete K\"ahler manifolds
satisfying Theorem 4.1.
\endproclaim

The following corollary is a simple consequence of the sharp
estimate. The interested reader might want to compare the
corollary with the example of \cite{Do}, for which there are more
sub-quadratic harmonic functions than the linear growth ones.

\proclaim{Corollary 4.1} Let $M$ be as in Theorem 4.1. Then, for
any $\e>0$,
$$
\dim({\Cal O}_{2-\e}(M)) \le m+1.
$$
\endproclaim

In the next result we combine the techniques in the proof of
Theorem 4.1 and Theorem 3.2 to obtain the dimension estimate for
the holomorphic sections of polynomial growth on line bundles with
controlled positive curvature. (We assume that $\Omega(x)$ is
continuous for the simplicity.)

\proclaim{Theorem 4.3} Let $M$ be a complete K\"ahler manifold
with nonnegative bisectional curvature. Let $(L, h)$ be a
holomorphic line bundle over $M$ such that $\|\Omega_{\abb}\|(x)$
satisfies (2.7). Assume that the curvature of $(L, h)$,
$\rho=\frac{\sqrt{-1}}{2\pi}\Omega_{\abb}\, dz_{\a}\wedge
d\bar{z}_\beta$ satisfies
$$
\nu(\Omega_{+}, x, \infty)=\limsup_{r\to \infty} \frac{1}{2m} \lf(
r^2\aint_{B_x(r)}\Omega_{+}(y)\, dv_y\ri) <\infty. \tag 4.20
$$
Then for any $s\in {\Cal O}_d(M, L)$
$$
mult_x([s])\le C_6\,\lf( d +\nu_\infty\ri). \tag 4.21
$$
Here we denote $\nu(\Omega_{+}, x, \infty)$ by $\nu_\infty$. In
particular, it implies
$$
\dim({\Cal O}_d(M, L))\le C_7(m) (d+\nu_\infty)^m. \tag 4.22
$$
If furthermore, we assume that $M$ has maximum volume growth, then
$$
mult_x([s])\le d+\nu_\infty. \tag 4.23
$$
In particular,
$$
\dim({\Cal O}_d(M, L))\le \dim_{\C}({\Cal
O}_{[d+\nu_\infty]}(\C^m)). \tag 4.24
$$
\endproclaim
\demo{Proof} We assume that ${\Cal O}_d(M, L)\ne \emptyset$.
Otherwise there is nothing to prove. As we observe in the remarks
after Corollary 3.3, it implies that
$$
\aint_{B_x(r)}\Omega_{-}(y)\, dv_y \le \frac{C_4}{r^2}.
$$
Combining with the assumption (4.20) we know that we have global
solution $v(x,t)$ to the Hermitian-Einstein heat equation (2.1)
(equivalently, (2.5)). We also consider the equation $\heatt
\tilde v(x,t) =0$ with initial data $\tilde v(x,0)=\D \log
\|s\|(x)$. By the approximation argument of Lemma 3.1, since both
$\Omega_{\abb}(x,t)$ and $\tilde v_{\abb}(x,t)$ satisfies (1.1),
and $\lf(\Omega_{\abb}+\tilde v_{\abb}\ri)(x,0)\ge 0$, by Lemma
2.1, we know that $\lf(\Omega_{\abb}+\tilde v_{\abb}\ri)(x,t)\ge
0$. Now we can apply Theorem 1.1 to $\lf(\Omega_{\abb}+\tilde
v_{\abb}\ri)(x,t)$. In particular we have that $(tw(x,t))_t\ge 0$.
Here $w=\Omega(x,t)+\tilde v_t(x,t)$. One can show that
$$
\lim_{t\to 0} t\, w(x,t) =\frac{1}{2} mult_x([s]) \tag 4.25
$$
 still holds since
$\Omega_{\abb}(x,0)=\Omega_{\abb}(x,0)$ does not contribute to the
Lelong number by the continuity assumption on $\Omega$. Similar to
the proof of (3.13) we have that
$$
\limsup_{t\to \infty} t\tilde v_t(x,t)\le C_1(m) \, d. \tag 3.13'
$$
This can be justified as follows. Since $\Omega(x,t)$ satisfies
the heat equation with $\Omega(x, 0)=\Omega(x)$, by Theorem 3.1 of
\cite{N1} we have that
$$
\Omega(x,t)\le C_7(m) \sup_{r\ge
\sqrt{t}}\lf(\aint_{B_x(r)}\Omega_{+}(y)\, dv_y\ri),
$$
which implies
$$
\limsup_{t\to \infty} t\Omega(x,t)\le C_7(m) \nu_\infty. \tag 4.26
$$
Here we have used (4.20). Combining (3.13') and (4.26) we get
(4.21). This finishes the proof for the general case.

For the case of maximum volume volume growth, we need to show that
$$
\lim_{t\to \infty} t\, w(x,t)\le \frac{1}{2}\lf(d+\nu_\infty\ri).
\tag 4.27
$$
Applying Lemma 4.2 to $\tv(x,t)$,  we then have that
$$
\limsup_{t\to\infty} t\tv_t(x,t)\le \frac{1}{2}d.
$$
The similar argument as the proof of  Lemma 4.2 also shows that
$$
\limsup_{t\to\infty}t\, \Omega(x,t)\le \frac{1}{2} \nu_\infty.
$$
Combining them together we have (4.27).
\enddemo

\proclaim{Remark} One can think $\nu(\Omega, x, \infty)$ is the
Lelong number of $\Omega$ at infinity. Similarly, one can define
the Lelong number of plurisubharmonic function at infinity. Then
Theorem 2.2 can be simply rephrased as that any nonnegative line
bundle has positive  Lelong number at infinity, if it is not flat,
and any nonconstant plurisubharmonic function  has positive Lelong
number at infinity, respectively.
\endproclaim

\input amstex
\documentstyle{amsppt}
\magnification=1200 \hsize=13.8cm \catcode`\@=11
\def\NoLogo{\let\logo@\empty}
\catcode`\@=\active \NoLogo
\def\heatt{\lf (\Delta-\frac{\p}{\p t}\ri)}

\def\heat{\lf(\frac{\p}{\p t}-\Delta\ri)}

\def \b {\beta}

\def\lf{\left}
\def\ri{\right}
\def\bbar{\bar \beta}
\def\a{\alpha}

\def\g{\gamma}
\def\e{\epsilon}
\def\p{\partial}
\def\delbar{{\bar\delta}}
\def\ddbar{\partial\bar\partial}
\def\dbar{\bar\partial}

\def\C{\Bbb C}
\def\N{\Bbb N}

\def\R{\Bbb R}
\def\tv{\tilde v}
\def\vp{\varphi}

\def\dbar{\bar\partial}

\def\abb{{\alpha\bar\beta}}
\def\gbd{{\gamma\bar\delta}}

\def \D {\Delta}
\def\aint{\frac{\ \ }{\ \ }{\hskip -0.4cm}\int}
\vsize=19.0 cm

\subheading{\S5 Applications}

\vskip .2cm

In this section we show some easy consequences of the multiplicity
estimates proved in Section 3. Further applications will be
treated later. We first define the holomorphic maps $\Phi_j$
inductively, for $j\in \N$, as follows.
$$
\Phi_1(x) =( f^1_1(x), \cdots, f^1_{k_1}(x))\in \C^{k_1}
$$
where $f^1_i$  form a basis for ${\Cal O}_1(M)/\C$. Suppose that
we have defined the map $\Phi_j:M\to \C^{k_j}$ as
$$
\Phi_j(x)=( f^j_1(x), \cdots, f^j_{k_j}(x)),
$$
where $k_j=\dim({\Cal O}_j (M)$ and $f^j_i$ form a basis of ${\Cal
O}_j(M)$. We define
$$
\Phi_{j+1}(x)=(f^{j+1}_1(x),\cdots, f^{j+1}_{k_j}(x),
f^{j+1}_{k_j+1}(x), \cdots, f^{j+1}_{k_{j+1}}(x))
$$
such that $f^{j+1}_i(x)=f^{j}_i(x)$ for all $i\le k_j$ and
$f^{j+1}_i$ with $i\ge k_j+1$ form a basis of ${\Cal
O}_{j+1}(M)/{\Cal O}_j(M)$. We define the Kodaira dimension $k(M)$
of $M$ as
$$k(M)=\max_{j\in \N}\{\max_{x\in M} (rank(\Phi_j(x)))\}. \tag 5.1
$$
In the case ${\Cal O}_P(M)=\C$ we define $k(M)=0$.

\proclaim{Proposition 5.1} Let $M$ be a complete K\"ahler manifold
with nonnegative bisectional curvature. Let ${\Cal M}(M)$ be the
quotient field generated by ${\Cal O}_P(M)$. Then
$$
deg_{tr}({\Cal M}(M)) =k(M). \tag 5.2
$$
\endproclaim
\demo{Proof} If $\max_{x\in M}(rank(\Phi_j))\ge k$, it is easy to
see that $deg_{tr}({\Cal M}(M))\ge k$ since there are at least $k$
holomorphic functions in ${\Cal O}_j(M)$ which are transcendental.
This implies that $deg_{tr}({\Cal M}(M))\ge k(M)$.

Now we show that $deg_{tr}({\Cal M}(M))\le k(M)$ too. By Theorem
3.1 we have that $ord_x(f)\le C_1(m) d$ for any $f\in {\Cal
O}_d(M)$. Assume that the conclusion is not true. Then we have
$F_1, \cdots, F_{k(M)+1}\in {\Cal M}(M)$ such that they are
transcendental over $\C$. We can assume that they have the form
$F_i= \frac{f_i}{f_0}$ with $f_j\in {\Cal O}_{d_0}(M)$. By
counting the monomial formed by $F_i$ we conclude that
$$\dim({\Cal O}_{d_0\, p}(M))\ge C_2(m) \lf(p\ri)^{k(M)+1} \tag 5.3 $$
for some positive constant $C_2(m)$.

 On the other hand, by the definition of $k(M)$ we know that
for $j\ge j_0$, there exists $y\in M$ such that $(rank
(\Phi_j)(y))=k(M)$. Then the  basis $\{f^j_i(x)\}$, $i=1,\cdots,
k_j$, are constant on $\Phi^{-1}(y)$. Namely $f^j_i(x)$ are
function of $k(M)$ free variables. Pick $x_0\in\Phi^{-1}(y)$. This
shows that it only takes dimension of $C_3(m)q^{k(M)}$ to get a
nontrivial function $f\in {\Cal O}_j(M)$ such that
$ord_{x_0}(f)\ge q$. Here $C_3(m)$ is a positive constant only
depends on $m$.

We can choose $d$ such that $d_0|d$ and
$$
C_2(m)\frac{d^{k(M)+1}}{d_0^{k(M)+1}}\ge C_3(m)(C_1(m)
d+1)^{k(M)}.
$$
From (5.3) this implies that there are $f\in {\Cal O}_d(M)$ such
that $ord_{x_0}(f)\ge C_1(m)d+1$. This is a contradiction.

\enddemo

\proclaim{Corollary 5.1} Let $F_1, \cdots, F_{k(M)}$ are the
transcendental elements in ${\Cal M}(M)$. Then the quotient field
${\Cal M}(M)$ is a finite algebraic extension of $\C(F_1,
\cdots,F_{k(M)})$. In particular, ${\Cal M}(M)$ is finitely
generated.
\endproclaim
\demo{Proof} The finiteness of the extension follows a similar
argument as Proposition 5.1. Once we know that the extension is
finite, the Primitive Element Theorem implies the finite
generation of ${\Cal M}(M)$ (cf. \cite{ZS, page 84}).
\enddemo

From the proof of Proposition 5.1 we also know that
\proclaim{Corollary 5.2} $$ \dim({\Cal O}_d(M))\le C_4(m)
d^{k(M)}. \tag 5.4
$$
\endproclaim

\proclaim{Remark} Result similar to Proposition 5.1 for the
compact manifolds, where the holomorphic functions in ${\Cal
O}_d(M)$ are replaced by the holomorphic sections of the power
$L^d$ of a fixed line bundle $L$,  was known as Serre-Siegel lemma
(cf. (6.5) in \cite{(D}). The dimension estimates as in Theorem
3.1 is much easier to prove in that case.
\endproclaim

It order to prove Theorem 0.3 we need the following well-known
result on the $L^2$-estimate of $\dbar$.

\proclaim{Theorem 5.1 [cf. \cite{A-V}, \cite{D}]} Let (E,g) be a
Hermitian line bundle with semi-positive curvature on complete
K\"ahler manifold $(M, h)$ of dimension m. Suppose $\varphi
:M\rightarrow [-\infty, 0]$ is a function of class $C^{\infty}$
outside a discrete subsets $s$ of $M$ and, near each point $p\in
S$, $\varphi(z)\ =\ A_p\log |z|^2$ where $A_p$ is a positive
constant and $z\ =\ (z_1,z_2, \dots, z_m)$ are local holomorphic
coordinates centered at $p$. Assume that $\Theta(E,
g\exp(-\varphi))\ =\ \Theta(E,g)\ + \ \p\bar{\p}\varphi \ge 0$ on
$M\setminus S$, and $\epsilon: M\rightarrow [0, 1]$ be a
continuous function such that $ \Theta(E,g)\ + \ \p\bar{\p}\varphi
\ge \epsilon \omega_{h}$ on $M\setminus S$. Then, for every
$C^{\infty}$ form $\theta$ of type (m,1) with values in $L$ on $M$
which satisfies
$$\bar{\p}\theta \ =\ 0 \ \text{and} \ \int_M \epsilon^{-1} |\theta|^2
e^{-\varphi}\ dv_h<\infty,$$ there exists a $C^{\infty}$ form
$\eta$ of type (m,0) with values in $L$ on $M$ such that
$$\bar{\p} \eta\ =\ \theta \ \text{and}\ \int_M |\eta|^2 e^{-\varphi}\ dv_h \le
\int_M \epsilon^{-1} |\theta|^2 e^{-\varphi}\ dv_h<\infty.$$
\endproclaim

\demo{Proof of Theorem 0.3} The proof uses a nice idea from
\cite{NR}. The estimate in Theorem 3.1 holds also the  key to make
the that argument work in this case. Let $\tilde M$ be a covering
space  of $M$. Denote $\pi$ the covering map. By the assumption
that $deg_{tr}({\Cal M}(M))=m$ we know $k(M)=m$, from Proposition
5.1. In particular, it implies that there exists a smooth
 plurisubharmonic function $\phi$, which can be constructed using
the holomorphic functions in ${\Cal O}_P(M)$ such that $\phi$ is
strictly plurisubharmonic at some point $p\in M$. Moreover, it
satisfies
$$
0\le \phi(x)\le C_5(M)\log( r(x)+2). \tag 5.5
$$
Here $r(x)$ is the distance function to a fixed point $o\in M$ and
$C_5(M)$ is a positive constant only depends on $M$. Now we choose
a small coordinate neighborhood $W$ near $p$ such that it is
evenly covered by $\pi$ such that $\sqrt{-1}\ddbar\phi>0$ in $W$.
Choose $p\in U\subset W$. Denote $U_{i}$ ($W_i$), $i=1, \cdots,l$
to be the disjointed pre-images of $U$ ($W$), and $p_i$ the
pre-images of $p$. Here $l$ is the number of the covering sheets,
which could be infinity. We use the coordinates $(z_1, \cdots,
z_m)$ in $W$, as well as in $U_i$ such that $p$ ($p_i$) is the
origin. Clearly, it does not hurt to assume that $W$ is inside the
ball $\{z|\, |z|\le 1\}$. We also define $\vp_p(x)=\rho(x) \log
|z|^2$, where $\rho(x)$ is a cut-off function with support inside
$W$, equal to $1$ in $U$. By changing the constant in (5.5), we
can make sure that $\phi+\vp_p$ is plurisubharmonic and
$$
\sqrt{-1}\ddbar \lf(\phi+\vp_p\ri)>0 \tag 5.6
$$
inside $U$. Now we lift $\phi$ to the  cover $\tilde M$ and denote
by $\tilde \phi$. Similarly $\tilde\vp_p$ is the lift of $\vp_p$.
Clearly (5.4) holds for $\tilde \phi$ and (5.5) holds for $\tilde
\phi+\tilde\vp_p$. Now we use Theorem 5.1 to show that for any
$d\in \N$, $1\le i\le l$, and $\a=(\a_1,\cdots, \a_m)$ with
$|\a|=d$, we can construct holomorphic function $f^i_{\a}(x)$ on
$\tilde M$ such that
$$
D^{\beta} f^i_{\a}\, (p_j)=\cases  1, \beta=\a, j=i,\\
0, \text{  otherwise} \endcases.
$$
To achieve this we first give such function locally in $U_i$,
which is trivial to do,  and then extend it by cutting-off to
whole $\tilde M$. Let's call it $\zeta$. Notice that we can
arrange $\zeta$ to be holomorphic in $U_i$. Now we let $\theta
=\dbar{\zeta}$. Apply Theorem 5.1 with $E=K_{\tilde M}^{-1}$ with
the metric $|\cdot|\exp(-\mu (\tilde \phi+\tilde\vp_p))$ (namely
$\vp=\tilde \vp_p$) with $\mu =m +d+3$. Theorem 5.1 then provides
$\eta$ such that $\dbar \eta =\theta$. Now $\zeta-\eta$ gives the
wanted holomorphic function, since by the choice of $\mu$, and the
finiteness of
$$
\int_M |\eta|^2\exp(-\mu(\tilde\phi+\tilde \vp_p))\, dv \tag 5.7
$$
implies  that $\eta$ vanishes at least up to order $d+1$ at $p_j$.
The finiteness of (5.7) also implies that
$$
\int_M |f_{\a}^i|^2\exp(-\mu(\tilde \phi+\tilde \vp_p))\, dv\le
B<\infty \tag 5.8
$$
for some positive constant $B$. Observe that $\vp_p\le 0$.
Combining with (5.5) we further have
$$
\int_{B_{\tilde o}(r)} |f^i_{\a}|^2\, dv\le B (r(x)+2)^{\mu C_5}.
\tag 5.9
$$
Using the mean value inequality of Li-Schoen we can conclude that
$f^i_\a \in {\Cal O}_{\frac{\mu C_5}{2}}(\tilde M)$. Noticing that
$f^i_\a$ are linearly independent, the  dimension of the space
spanned by them is bounded from below by $C_6(m)d^m l$. But from
Theorem 3.1 we also know that $$\dim({\Cal O}_{\frac{\mu
C_5}{2}}(\tilde M))\le C_7(m)\lf(\frac{\mu C_5}{2}\ri)^m.$$
Plugging $\mu =d+m+3$ we have that
$$
C_6(m)d^m l\le C_7(m)\lf(\frac{ C_5(d+m+3)}{2}\ri)^m
$$
which implies $l\le \frac{C_7(m)}{C_6(m)}\lf(\frac{C_5}{2}\ri)^m$,
by letting $d\to \infty$.

\enddemo

\proclaim{Corollary 5.3} Let $M$ be a complete K\"ahler manifolds
with nonnegative bisectional curvature. Assume that the Ricci is
positive some where and the scalar curvature ${\Cal R}(x)$
satisfies (2.7) and (3.24). Namely
$$
\sup_{r\ge 0}\lf(\exp(-ar^2)\aint_{B_o(r)}{\Cal R}(y)\,
dy\ri)<\infty
$$
and
$$
\aint_{B_o(r)}{\Cal R}(y)\, dv_y \le \frac{C}{r^2}
$$
for some positive constants $a$ and $C$. Then $M$ has finite
fundamental group.
\endproclaim
\demo{Proof} The result follows from Theorem 0.3 and Corollary 6.2
of \cite{NT2}.
\enddemo

\enddocument